\documentclass[12pt]{article}
\RequirePackage{amsthm, amsmath, amsfonts, amssymb, graphicx}
\RequirePackage[numbers,square]{natbib}
\usepackage{dsfont}
\usepackage{setspace}
\newcommand{\bneqn}{\vspace{-0.25cm}\begin{eqnarray}}
\newcommand{\eneqn}{\end{eqnarray}}
\allowdisplaybreaks

\usepackage[utf8]{inputenc}
\usepackage{hyperref,bookmark}
\hypersetup{
  colorlinks=true,
  linkcolor=red,
  citecolor=blue,
  filecolor=blue,
  urlcolor=blue,
}

\newtheorem{theorem}{Theorem}
\newtheorem{lemma}[theorem]{Lemma}
\newtheorem{remark}{Remark}
\newtheorem{proposition}{Proposition}
\newtheorem{corollary}{Corollary}
\newcommand{\abss}[1]{\left|#1\right|}
\newcommand{\tabincell}[2]{\begin{tabular}{@{}#1@{}}#2\end{tabular}}

\advance \topmargin by -\headheight
\advance \topmargin by -\headsep
\advance \topmargin .2in
\title{The rencontre problem}

\title{The rencontre problem}
\author{F. Thomas Bruss\footnote{Département de Math\'ematique,
Université Libre de Bruxelles}, Philip A. Ernst\footnote{Department of Statistics, Rice University}, and Dongzhou Huang\footnote{Department of Statistics, Rice University}  \\\\
\textbf{In Memory of Larry Shepp}}

\begin{document}
\maketitle

\begin{abstract}
\noindent Let $\left\{X^{1}_k\right\}_{k=1}^{\infty}, \left\{X^{2}_k\right\}_{k=1}^{\infty}, \cdots, \left\{X^{d}_k\right\}_{k=1}^{\infty}$ be $d$ independent sequences of Bernoulli random variables with success-parameters $p_1, p_2, \cdots, p_d$ respectively, where $d \geq 2$ is a positive integer, and $ 0<p_j<1$ for all $j=1,2,\cdots,d.$ Let
\begin{equation*}
 S^{j}(n) = \sum_{i=1}^{n} X^{j}_{i} = X^{j}_{1} + X^{j}_{2} + \cdots + X^{j}_{n}, \quad n =1,2 , \cdots.
\end{equation*}
We declare a ``rencontre'' at time $n$, or, equivalently, say that $n$ is a ``rencontre-time,'' if
\begin{equation*}
S^{1}(n) = S^{2}(n) = \cdots = S^{d}(n).
\end{equation*}
We motivate and study the distribution of the \textit{first} (provided it is finite) rencontre time.
\medskip



\end{abstract}

\section{Introduction}

Consider $\left\{X^{1}_k\right\}_{k=1}^{\infty}, \left\{X^{2}_k\right\}_{k=1}^{\infty}, \cdots, \left\{X^{d}_k\right\}_{k=1}^{\infty}$ to be $d$ independent sequences of Bernoulli random variables with success-parameters $p_1, p_2, \cdots, p_d$ respectively, where $d \geq 2$ is a positive integer, and $0<p_j<1$ for all $j=1, 2, \cdots, d.$ Let
\begin{equation*}
 S^{j}(n) = \sum_{i=1}^{n} X^{j}_{i} = X^{j}_{1} + X^{j}_{2} + \cdots + X^{j}_{n}, \quad n =1,2 , \cdots.
\end{equation*}
We declare a ``rencontre'' at time $n$, or equivalently, say that $n$ is a ``rencontre-time,'' if $n\ge1$ and
\begin{equation*}
S^{1}(n) = S^{2}(n) = \cdots = S^{d}(n).
\end{equation*}
In plain English, the event that there is a rencontre at time $n \geq 1$ is exactly the event $\{ S^{1}(n) = S^{2}(n) = \cdots = S^{d}(n) \}$.
The first rencontre-time is given as
\begin{equation*}
J^d:= J^d(p_1, p_2, \cdots, p_d) = \inf\{{n \in \{1,2,\cdots\} : n \text{ a rencontre-time}}\},
\end{equation*}
that is, $J^d$ is the first time the random walk $\left\{\left(S^1_n, \cdots, S^d_n \right) : n\geq1 \right\}$ intersects with the line $\left\{(x_1,\cdots,x_d) : x_1  =\cdots = x_d  \right\}$. Further, let $q_j = 1- p_j, j=1,2,\cdots, d$. In order to exclude trivialities, or evident remarks about possible reduction of dimension $d$, we shall suppose that all parameters $p_1, p_2, \cdots, p_d$ are strictly between $0$ and $1.$ The present work studies the distribution of the first rencontre time $J^d$ (provided such a time exists).

\smallskip

\indent  We shall see that the case $d=2$ is special in the sense that, when $p_1=p_2, P(J^2<\infty)=1$ and, for all values of $p_1$ and $p_2$,
$E(J^2)=\infty.$  By a simple projection argument we may conclude without any further calculations that $E(J^d)=\infty$ for $d\ge 2.$ Indeed, in order to have a rencontre at some time $t$ it is necessary to have a rencontre in all $d \choose 2$ different pairs of the defined Bernoulli processes, so that $$E(J^d) \ge \max\{E(J^2_{k,\ell}): 1\le k<\ell\le d\}=\infty,$$
 where $J^2_{k,\ell}$ denotes the corresponding first rencontre time for the $k$th and $\ell$th subprocess. This is why our main interest shall be on the distribution.
We also remark that although the general problem can be converted to the problem of first intersection to the origin of $(d-1)$-dimensional random walks by considering $\widetilde{S}_{n} = \left(S^1_n - S^d_n, \cdots, S^{d-1}_n- S^{d}_{n} \right)$, this formulation proves more unwieldy. \\
\indent  The literature most closely related to this problem studies the number of intersections of $n$ independent simple random walks. For two processes $\{S_n\}$ and $\{T_n\}$, references \cite{Pemantle, Lawler, Szekely} consider the cardinality of the set $\{k \in \mathbb{N}: k=S_n=T_m\,\, \text{for any} \,\, m,n\}$. Our paper departs from these previous works in that we are only interested in the \textit{first time} of intersection.

\indent We now offer two practical motivations for the problem we consider:
\begin{enumerate}
\item Consider $d$ independent sequences of Bernoulli random variables with success-parameter $p_1, p_2, \cdots, p_d$ respectively, for $d \geq 2$ a positive integer. Suppose that the sequences model strands of genes and that a zero is assigned to a gene which is not activated and a one is assigned to a gene which is activated. We may be interested in the \textit{first time} when the number of activated genes coincides across these sequences.
\item Suppose that two players, A and B, play a sequence of independent games with each other. Let $p_A$ be the win probability for player A in any given game, $p_B$ be the win probability for player B in any given game, each independently of each other. Let $S_A(n)$ and $S_B(n)$ be the respective scores of players A and B after $n$ rounds. Now suppose that both players A and B can quit the game without cost at a rencontre-time, that is at  the time $t$ such that $S_A(t)=S_B(t)$. Further suppose that the current loser at time $t'$ would have to pay $\abss{S_A(t')-S_B(t')}$. It now becomes of interest to know the distribution of the waiting time until the next rencontre-time.
\end{enumerate}

\indent The remainder of this manuscript is organized as follows. Section \ref{sec2} derives and discusses the distribution of $J^d$ (the first rencontre-time). In Section \ref{sec3}, we introduce the probability generating function of $J^d$ and present a link between the latter and the generating function of probabilities of having a rencontre at any given time.  In Section \ref{secexp}, we derive an explicit form of probability generating function of $J^d$ and use characteristic functions in order to provide an expression for $P(J^d = \infty)$. In Section \ref{sec5}, we give an alternative proof (Theorem \ref{theorem:4}) that the expectation of $J^d$ is infinite for $d \geq 3$. This is clear from our preceding result for $d=2$ and the projection argument given above. However this alternative proof of Theorem 5 offers
a clear benefit providing estimates which are useful for
estimating the conditional expectations $E(J^d | J^d < \infty)$ and $E(J^d|b < J^d<\infty)$ for some upper bound $b$. We pursue this task in Section \ref{sec6}.
\section{Distribution of the first rencontre-time} \label{sec2}
We say that a rencontre happens at time $n$ in state $k$ if
\begin{equation*}
 \left(S^{1}(n), S^{2}(n), \cdots, S^{d}(n)\right) = (k,k,\cdots, k).
\end{equation*}
Note that this definition implies that $k \in \{0, 1,2, \cdots,n\}$. Since the i.i.d. random walks are independent of each other, we have that
\begin{equation*}
P(\text{rencontre at time $n$ in state $k$} )
= \prod_{j=1}^{d} P\left(S^{j}(n)=k\right)
= \prod_{j=1}^{d} {n \choose k} p_j^k q_j^{n-k}.
\end{equation*}
Let $R^d_n, n=1,2\cdots$ denote the event that a rencontre happens at time $n$ for these $d$ random walks. Thus,\,$R^d_n$ may be written as union of disjoint events as
\begin{equation*}
R^d_{n} = \bigcup_{k=0}^{n}\{ \text{rencontre at time $n$ in state $k$} \}.
\end{equation*}
It then follows that
\begin{equation}
P(R^d_n) = \sum_{k=0}^{n} P(\text{rencontre at time $n$ in state $k$})
= \sum_{k=0}^{n} \prod_{j=1}^{d} {n \choose k} p_j^k q_j^{n-k}.  \label{formforR}
\end{equation}
We now proceed with Theorem \ref{theorem:1}, which indeed is an instance of ``first-occurrence decomposition'' in Feller's theory of recurrent events (\cite{Feller}).
\begin{theorem} \label{theorem:1}
For $n\in \mathds{N}_{+}$, we have
\begin{equation}
P(J^d=n) = \sum_{s=1}^{n} (-1)^{s-1} \sum_{j_1+\cdots+j_s=n} P(R^d_{j_1})\cdots P(R^d_{j_s}). \label{main_equation}
\end{equation}
\end{theorem}
\begin{proof}
\begin{eqnarray*}
\{J^d=n\} &=& \{\text{no rencontre up to time $n-1$, rencontre at time $n$}\} \\
 &=& R^d_n \backslash \bigcup_{s=1}^{n-1} R^d_s
 = R^d_n \backslash \bigcup_{s=1}^{n-1} \left(R^d_s \cap R^d_n \right).
\end{eqnarray*}
The probability of the event $J^d =n$ is
\begin{eqnarray}
P(J^d=n)&=& P\bigg(R^d_n \backslash \bigcup_{s=1}^{n-1} \left(R^d_s \cap R^d_n \right)\bigg) \notag \\
    &=& P(R^d_n) - P\left(\bigcup_{s=1}^{n-1} \left(R^d_s \cap R^d_n \right)\right).  \label{equation:1}
\end{eqnarray}
By inclusion-exclusion, we have
\begin{eqnarray}
&&P\left(\bigcup_{s=1}^{n-1} \left(R^d_s \cap R^d_n \right)\right) \notag \\
&=& \sum_{s=1}^{n-1} (-1)^{s-1} \sum_{1\leq j_1<\cdots < j_s\leq n-1} P\left( \left(R^d_{j_1}\cap R^d_{n}\Big) \cap \cdots \cap \Big(R^d_{j_s}\cap R^d_{n}\right) \right)\notag \\
&=&\sum_{s=1}^{n-1} (-1)^{s-1} \sum_{1 \leq j_1<\cdots < j_s\leq n-1} P\left( R^d_{j_1} \cap \cdots \cap R^d_{j_s}\cap R^d_{n} \right). \label{equation:2}
\end{eqnarray}
We shall use recursive arguments to simplify the probability of intersection of events in \eqref{equation:2}. For example, for $j_1 < j_2$,
\begin{equation*}
P\left(R^d_{j_1}\cap R^d_{j_2}\right) = P\left( R^d_{j_1} \right)P\left( R^d_{j_2-j_1} \right).
\end{equation*}
Knowledge of a rencontre at time $j_1$ allows the $d$ processes to be in the same state (and, for simplicity, we may consider them all as starting again from $(0,0,\cdots,0)$). By induction, the terms in \eqref{equation:2} split into the corresponding product
\begin{equation} \label{anewone}
P\left( R^d_{j_1} \cap \cdots \cap R^d_{j_s}\cap R^d_{n} \right)
=P\left( R^d_{j_1}\right) P\left(R^d_{j_2-j_1}\right)\cdots P\left(R^d_{j_s-j_{s-1}}\right)  P\left(R^d_{n-j_s}\right).
\end{equation}
Plugging (\ref{anewone}) into \eqref{equation:2} gives
\begin{eqnarray}
&&P\left(\bigcup_{s=1}^{n-1} \left(R^d_s \cap R^d_n \right)\right) \notag \\
&=& \sum_{s=1}^{n-1} (-1)^{s-1} \sum_{1 \leq j_1<\cdots < j_s\leq n-1} P( R^d_{j_1}) P(R^d_{j_2-j_1})\cdots P(R^d_{n-j_s}) \label{equation:3}.
\end{eqnarray}
Let $l_u = j_u - j_{u-1}$, $u\leq s$ and $l_{s+1}= n - j_s$, where by convention $j_0=0$. The right-hand side of equation (\ref{equation:3}) simplifies to
\begin{equation}
\sum_{s=1}^{n-1} (-1)^{s-1} \sum_{l_1+\cdots+l_{s+1}=n} P( R^d_{l_1}) P(R^d_{l_2})\cdots P(R^d_{l_{s+1}})\label{equation:4}.
\end{equation}
We now perform a change of variables $\tilde{s}= s+1$. The right-hand side now simplifies to
\begin{equation} \sum_{s=2}^{n} (-1)^{s} \sum_{l_1+\cdots+l_{s}=n} P( R^d_{l_1}) P(R^d_{l_2})\cdots P(R^d_{l_{s}}). \label{equation:5}
\end{equation}
Combining \eqref{equation:1} and \eqref{equation:5} completes the proof.
\end{proof}

\section{Probability generating function of $J^d$} \label{sec3}
Theorem \ref{theorem:1} provides an expression for $P(J^d=n)$ but does not allow us to compute $P(J^d=\infty)$ (i.e., the probability of no rencontre). 
We hence turn to generating functions. Let us define
\begin{equation} 
\phi_{d}(x) := \phi_{d}(x;p_1,\cdots,p_d)= \sum_{n=1}^{\infty} P\left(J^d=n\right) x^n, \label{equation:6}
\end{equation}
and
\begin{equation}
\varphi_{d}(x) := \varphi_{d}(x;p_1,\cdots,p_d) = \sum_{n=1}^{\infty} P\left( R^{d}_n \right) x^n. \label{equation:7}
\end{equation}
Note that since $\sum_{n=1}^{\infty} P\left(J^d=n\right) \leq 1$, the power series in \eqref{equation:6} converges if $x \in [0,1]$. For $P\left( R^{d}_n \right) \leq 1$, the power series in \eqref{equation:7} converges if $x \in [0,1) $. Recursive arguments enables us to show that $\phi_{d}(x)$ is related to $\varphi_{d}(x)$ as follows:
\begin{lemma} \label{proposition:1}
For $x \in [0,1)$, we have
\begin{equation*}
 1 - \phi_{d}(x) = \frac{1}{1+ \varphi_{d}(x) }.
\end{equation*}
\end{lemma}
\begin{proof}
This Lemma is an instance of the ``Feller relation'' and is proven in Theorem $1$ in Chapter $13.3$ of Feller (\cite{Feller}). Note that Feller's $F$ is our $\phi_{d}$ and Feller's $U$ is our $1+ \varphi_{d}$.
\end{proof}

\section{An expression for $P(J^d = \infty)$} \label{secexp}
Note that the coefficients in the power series in \eqref{equation:6} are non-negative. By Abel's theorem for power series, we have
\begin{equation*}
\sum_{n=1}^{\infty}P\left(J^d =n\right) = \lim_{x \rightarrow 1-} \sum_{n=1}^{\infty} P\left(J^d =n\right) x^n
= \lim_{x \rightarrow 1-} \phi_{d}(x),
\end{equation*}
since by definition $\sum_{n=1}^{\infty} P\left(J^d =n\right) \leq 1$. Similarly,
\begin{equation}
\sum_{n=1}^{\infty} P\left( R^d_n \right) = \lim_{x \rightarrow 1-} \sum_{n=1}^{\infty} P\left( R^d_n \right) x^n
= \lim_{x \rightarrow 1-} \varphi_{d}(x) = \varphi_{d}(1-). \label{varphileftlimt}
\end{equation}
Applying Lemma \ref{proposition:1} gives
\begin{eqnarray}
P(J^d = \infty) &=& 1- \sum_{n=1}^{\infty} P(J^d = n) = 1 - \lim_{x \rightarrow 1-} \phi_{d}(x) \notag \\
&=& \lim_{x \rightarrow 1-} \frac{1}{1+ \varphi_{d}(x)} = \frac{1}{1+\varphi_{d}(1-)} .  \label{equation:12}
\end{eqnarray}
This allows us to convert the problem of calculating $P(J^d = \infty)$ into the problem of calculating $1+\varphi_{d}(1-)$.
\subsection{Characteristic function representation}
\indent We shall now use characteristic functions to give an expression for $1+\varphi_{d}(x)$. Let $\underline{\theta}^{d}$ be the vector $\left( \theta_1, \cdots, \theta_d\right)$ and let $\underline{S^d_n}$ the vector $\left( S^1(n), \cdots, S^d(n) \right)$. For simplicity, we will write $\underline{\theta}^{d}$ as $\underline{\theta}$ and $\underline{S^d_n}$ as $ \underline{S_n}$. Let
\begin{equation*}
\psi_{d} \left(\underline{\theta}\right) := \psi_{d} \left(\underline{\theta}; p_1, \cdots, p_d\right)
\end{equation*}
be the characteristic function of $\underline{S_1}$ (i.e. $\left( X^1_1, \cdots, X^d_1 \right)$). Direct calculation gives
\begin{eqnarray*}
\psi_{d} \left(\underline{\theta}\right) &=& E \left(e^{i\, \underline{\theta} \left(\underline{S_1}\right)^{T} } \right) = E\left( e^{\sum_{j=1}^{d} i\,\theta_{j} X^{j}_1  } \right) \\
&=& \prod_{j=1}^{d} E\left( e^{i\,\theta_{j} X^{j}_1 } \right) = \prod_{j=1}^{d} \left( p_j \, e^{i\,\theta_j}+q_j \right).
\end{eqnarray*}
Let
\begin{equation*}
\psi_{d,n} \left( \underline{\theta} \right) :=  \psi_{d,n} \left( \underline{\theta};p_1, \cdots, p_d \right)
\end{equation*}
be the characteristic function of $\underline{S_n}$.
Since $\left\{X^{1}_k\right\}_{k=1}^{\infty}, \left\{X^{2}_k\right\}_{k=1}^{\infty}, \cdots, \left\{X^{d}_k\right\}_{k=1}^{\infty}$ are independent, and $\{X_k^j\}$ is a sequence of i.i.d. Bernoulli random variables, we have
\begin{eqnarray*}
\psi_{d,n}\left(\underline{\theta}\right)
= E \left(e^{i\, \underline{\theta} \left(\underline{S_n} \right)^{T}} \right)
= \left(E \left(e^{i\, \underline{\theta} \left(\underline{S_1}\right)^{T} } \right)\right)^n
= (\psi_{d} \left(\underline{\theta}\right))^n.
\end{eqnarray*}
The inversion formula for the characteristic function $\psi_{d,n}\left(\underline{\theta}\right)$ is
\begin{eqnarray*}
P\left(\underline{S_n}= (x_1,\cdots,x_d) \right)
&=& \frac{1}{(2\pi)^d}\int\cdots\int_{[-\pi,\pi]^{d}} e^{-i\,(x_1,\cdots,x_d)\,\left(\underline{\theta}\right)^{T}}
\cdot\psi_{d,n}\left(\underline{\theta}\right) d\underline{\theta} \\
&=& \frac{1}{(2\pi)^d}\int\cdots\int_{[-\pi,\pi]^{d}} e^{-i\,\sum_{j=1}^{d} x_j\,\theta_j }
\cdot\psi_{d,n}\left(\underline{\theta}\right) d\underline{\theta}.
\end{eqnarray*}
This formula gives us an additional expression for the probability of a rencontre at time $n$, i.e.
\begin{eqnarray*}
&& P(R^d_n) = \sum_{k=0}^{n} P(\underline{S_n}= (k,\cdots,k)) \\
&=& \sum_{k=0}^{n} \frac{1}{(2\pi)^d}\int\cdots\int_{[-\pi,\pi]^{d}}\  e^{-i\,\sum_{j=1}^{d} k\,\theta_j }
\cdot\psi_{d,n}\left(\underline{\theta}\right) d\underline{\theta} \\
&=&  \frac{1}{(2\pi)^d}\int\cdots\int_{[-\pi,\pi]^{d}} \  \sum_{k=0}^{n} e^{-ik\,\sum_{j=1}^{d}\,\theta_j }
\cdot\left(\psi_{d}\left(\underline{\theta}\right)\right)^n d\underline{\theta}.
\end{eqnarray*}
Note that $|p_j e^{i\theta_j} + q_j|\leq p_j|e^{i\theta_j}| + q_j=1$, and thus $|\psi_{d}\left(\underline{\theta}\right)|\leq 1$. For $x \in [0,1)$, by Dominated Convergence, we have
\begin{eqnarray}
 &&1+ \varphi_{d}(x) = 1+ \sum_{n=1}^{\infty} P(R^d_n) x^n \notag \\
 &=&1+ \sum_{n=1}^{\infty} x^n\, \frac{1}{(2\pi)^d}\int\cdots\int_{[-\pi,\pi]^{d}} \  \sum_{k=0}^{n} e^{-ik\,\sum_{j=1}^{d}\,\theta_j }
\cdot\left(\psi_{d}\left(\underline{\theta}\right)\right)^n \,d\underline{\theta} \notag \\
&=& 1+ \sum_{n=1}^{\infty} \, \frac{1}{(2\pi)^d}\int\cdots\int_{[-\pi,\pi]^{d}} \  \sum_{k=0}^{n} e^{-ik\,\sum_{j=1}^{d}\,\theta_j }
\cdot\left(x\,\psi_{d}\left(\underline{\theta}\right)\right)^n \,d\underline{\theta}  \notag  \\
&=& \sum_{n=0}^{\infty} \, \frac{1}{(2\pi)^d}\int\cdots\int_{[-\pi,\pi]^{d}} \  \sum_{k=0}^{n} e^{-ik\,\sum_{j=1}^{d}\,\theta_j }
\cdot\left(x\,\psi_{d}\left(\underline{\theta}\right)\right)^n \,d\underline{\theta} \notag \\
&=& \frac{1}{(2\pi)^d}\int\cdots\int_{[-\pi,\pi]^{d}} \ \sum_{n=0}^{\infty} \sum_{k=0}^{n} e^{-ik\,\sum_{j=1}^{d}\,\theta_j }
\cdot\left(x\,\psi_{d}\left(\underline{\theta}\right)\right)^n \, d\underline{\theta}  \notag \\
&=& \frac{1}{(2\pi)^d}\int\cdots\int_{[-\pi,\pi]^{d}} \  \sum_{k=0}^{\infty} \sum_{n=k}^{\infty} e^{-ik\,\sum_{j=1}^{d}\,\theta_j }
\cdot\left(x\,\psi_{d}\left(\underline{\theta}\right)\right)^n \, d\underline{\theta} \notag \\
&=& \frac{1}{(2\pi)^d}\int\cdots\int_{[-\pi,\pi]^{d}} \  \sum_{k=0}^{\infty} \,e^{-ik\,\sum_{j=1}^{d}\,\theta_j }
\cdot \frac{\left(x\,\psi_{d}\left(\underline{\theta}\right)\right)^k}{1-x\,\psi_{d}\left(\underline{\theta}\right)} \,d\underline{\theta} \notag \\
&=& \frac{1}{(2\pi)^d}\int\cdots\int_{[-\pi,\pi]^{d}} \   \frac{1}{1-x\,\psi_{d}\left(\underline{\theta}\right)}
\, \sum_{k=0}^{\infty} \left(  x\,\psi_{d}\left(\underline{\theta}\right) \, e^{-i \sum_{j=1}^{d} \theta_j}  \right)^k
\,d\underline{\theta} \notag \\
&=& \frac{1}{(2\pi)^d}\int\cdots\int_{[-\pi,\pi]^{d}} \   \frac{1}{\left(1-x\,\psi_{d}\left(\underline{\theta}\right)\right) \left(1- x\,\psi_{d}\left(\underline{\theta}\right) \, e^{-i \sum_{j=1}^{d} \theta_j} \right)}
\,d\underline{\theta}. \label{equation:13}
\end{eqnarray}
Together with \eqref{equation:12}, the above allows us to give an expression for $P(J^d=\infty)$ as follows:
\begin{equation}
P\left( J^d = \infty\right) = \lim_{x \rightarrow 1-}
\left(\frac{1}{(2\pi)^d}\int\cdots\int_{[-\pi,\pi]^{d}} \   \frac{1}{\left(1-x\,\psi_{d}\left(\underline{\theta}\right)\right) \left(1- x\,\psi_{d}\left(\underline{\theta}\right) \, e^{-i \sum_{j=1}^{d} \theta_j} \right)}
\,d\underline{\theta}\right)^{-1}. \label{equation:14}
\end{equation}
In Appendix A, we show in the case $d=2$, the function $1+\varphi_2(x)$ can be calculated explicitly as
\begin{equation}
1+\varphi_2(x) = \frac{1}{\sqrt{1- 2x(p_1 p_2 +q_1 q_2)+ x^2(p_1p_2-q_1q_2)^2}}, \label{expforvarphi2}
\end{equation}
and thus
\begin{equation}
\phi_{2} (x) = 1-  \sqrt{1- 2x(p_1 p_2 +q_1 q_2)+ x^2(p_1p_2-q_1q_2)^2}. \label{equation:15}
\end{equation}
In the case $d=2$, our model can be converted to one-dimensional random walk with a stay (i.e. the values of increment are $-1, 0 ,1$) by letting $\widetilde{S}_{n} =S^{1}_n - S^2_n=  \sum_{i=1}^{n} \left(X^1_{i}-X^2_i\right) $. Then the problem of a first rencontre is equivalent to problem of first return to $0$. The authors of \cite{Saroj1} considered the one-dimensional random walk with a stay in the presence of partially reflecting barriers $a$ and $-b$. Indeed, \eqref{expforvarphi2} is a special case of the results of \cite{Saroj1}. \\
\indent Recall from \eqref{equation:6} that $\phi_{2}(1)=\sum_{n=1}^\infty P(J^2=n)=P(J^2<\infty)$ so that $P(J^2=\infty)=1-\phi_{2}(1)$. It is now straightforward to check that $1-\phi_{2}(1)$ gives the following form in \eqref{equation:15}: $\sqrt{(p_1-p_2)^2}=\abss{p_1-p_2}$. We thus obtain Theorem \ref{theorem:2} below.
\begin{theorem} \label{theorem:2}
In the case $d=2$, i.e. two i.i.d. random walks which are independent of each other, the probability of no rencontre is $P\left(J^{2}=\infty\right)=|p_1 - p_2|$. For all $p_1$ and $p_2$, the expectation of $J^2$ is $E\left(J^2\right)= \infty$.
\end{theorem}

\section{Some estimation results} \label{sec5}
In equation \eqref{equation:14} of Section \ref{secexp}, we gave an expression for $P\left(J^d = \infty \right)$. However, the integral cannot be calculated explicitly. This makes it difficult to answer questions such as whether $P\left( J^d = \infty \right)$ (the probability of no rencontre) is zero or non-zero. The present section develops tools to answer this question. Note that by \eqref{equation:12}, we have
\begin{equation*}
 P\left( J^d = \infty \right)= \frac{1}{1+ \varphi_{d}(1-)},
\end{equation*}
which implies that $P\left( J^d = \infty \right) =0 $ if and only if $\varphi_{d}(1-)= \infty$. Combining equations \eqref{formforR} and \eqref{equation:7} gives
\begin{eqnarray*}
\varphi_{d}(x) &=& \sum_{n=1}^{\infty} P\left( R^d_n \right)\,x^n
= \sum_{n=1}^{\infty} x^n \, \sum_{k=0}^{n} \prod_{j=1}^{d} {n \choose k} p_{j}^{k} q_{j}^{n-k}  \\
&=& \sum_{n=1}^{\infty} x^n \, \left(\prod_{j=1}^{d} q_j \right)^n \sum_{k=0}^{n} {n \choose k}^d \left(\prod_{j=1}^{d}p_j q_{j}^{-1} \right)^k.
\end{eqnarray*}
Let $Q_d$ denote $\prod_{j=1}^{d} q_j$ and $P_d$ denote $\prod_{j=1}^{d}p_j q_{j}^{-1}$. For ease of notation, we will write $Q_d$ as $Q$ and $P_d$ as $P$. Then
\begin{equation}
\varphi_{d}(x) = \sum_{n=1}^{\infty} x^n Q^n \sum_{k=0}^{n} {n \choose k}^d P^k \label{eqforphi}.
\end{equation}
By Abel's theorem for power series,
\begin{equation}
\varphi_{d}(1-) = \sum_{n=1}^{\infty} Q^n \sum_{k=0}^{n} {n \choose k}^d P^k. \label{equation:16}
\end{equation}
In order to study the finiteness of $\varphi_{d}(1-)$, we need to estimate $ \sum_{k=0}^{n} {n \choose k}^d P^k $. In the sequel, we will give upper bounds and lower bounds for $ \sum_{k=0}^{n} {n \choose k}^d P^k $ for sufficiently large $n$. To find such bounds, we must provide a few propositions. The value of $\alpha$ in the forthcoming propositions is always assumed positive.

\begin{proposition} \label{proposition4}
\begin{spacing}{1.3}
Viewing ${n \choose k} \alpha^{k}$ as a function of $k$, $k \in \{0,1,\cdots,n\}$, then ${n \choose k} \alpha^{k}$ is non-decreasing if $k \in \big\{0,1, \cdots, \big[ \frac{\alpha(n+1)}{\alpha +1}\big]\big\}$ and non-increasing if $ k \in \big\{\big[ \frac{\alpha(n+1)}{\alpha +1}\big], \big[ \frac{\alpha(n+1)}{\alpha +1}\big]+1, \cdots,n \big\}$, where $[x]$ is the the greatest integer less than or equal to $x$. As a result, when $k=\big[ \frac{\alpha(n+1)}{\alpha +1}\big]$, ${n \choose k} \alpha^{k}$ obtains its maximum, i.e.
\end{spacing}
\begin{equation*}
 {n \choose k} \alpha^{k} \leq {n \choose \big[ \frac{\alpha(n+1)}{\alpha +1}\big]}\alpha^{\big[ \frac{\alpha(n+1)}{\alpha +1}\big]},
  \quad k \in \{0,1,\cdots,n\}.
\end{equation*}
\end{proposition}
\begin{proof}
\begin{equation}
\frac{{n \choose k+1} \alpha^{k+1}}{{n \choose k} \alpha^{k}}
= \frac{n-k}{k+1} \, \alpha, \label{equation:19}
\end{equation}
which is a decreasing function of $k$.
We set the right-hand of \eqref{equation:19} $\geq$ 1 and obtain
\begin{equation*}
 k \leq \frac{\alpha(n+1)}{\alpha +1} -1.
\end{equation*}
This concludes the proof.
\end{proof}
\begin{proposition} \label{proposition5}
For sufficiently large $n$, we have
\begin{equation}
 {n \choose \big[ \frac{\alpha(n+1)}{\alpha +1}\big]}\alpha^{\big[ \frac{\alpha(n+1)}{\alpha +1}\big]}
 = \left( \frac{\alpha+1}{\sqrt{2 \pi \, \alpha}} + o(1) \right) n^{-\frac{1}{2}} (1+\alpha)^{n}.
\end{equation}
\end{proposition}
\begin{proof}
Let $\beta$ denote $\big[ \frac{\alpha(n+1)}{\alpha +1} \big] $. It then follows that, for sufficient large $n$,
\begin{eqnarray} \label{eq25}
 \beta &=& \left( \frac{\alpha}{\alpha +1} +o(1) \right) n,  \nonumber \\
  n-\beta &=& \left( \frac{1}{\alpha +1} +o(1) \right) n ,
\end{eqnarray}
By Stirling's formula, we have
\begin{eqnarray}
&&{n \choose \big[ \frac{\alpha(n+1)}{\alpha +1}\big]}\alpha^{\big[ \frac{\alpha(n+1)}{\alpha +1}\big]}
= {n \choose \beta} \alpha^{\beta} = \frac{n!}{\beta! (n-\beta)!} \alpha^{\beta} \notag \\ [2mm]
&\sim& \frac{\sqrt{2 \pi n} \Big(\frac{n}{e}\Big)^n}{\sqrt{2 \pi \beta} \Big(\frac{\beta}{e}\Big)^\beta  \sqrt{2 \pi (n-\beta)} \Big(\frac{n-\beta}{e}\Big)^{n-\beta}} \, \alpha^{\beta} \notag  \\ [1mm]
&\sim& \frac{1}{\sqrt{2\pi}} \sqrt{\frac{n}{\beta(n-\beta)}}\left(\frac{n}{\beta}\right)^{\beta} \left(\frac{n}{n-\beta}\right)^{n-\beta} \, \alpha^{\beta} \notag \\ [1mm]
&\sim&  \frac{\alpha+1}{\sqrt{2 \pi \, \alpha}} n^{-\frac{1}{2}} \left(\frac{n}{\beta}\right)^{\beta} \left(\frac{n}{n-\beta}\right)^{n-\beta} \, \alpha^{\beta} \quad (\text{by}\, (\ref{eq25})) \notag \\ [1mm]
&\sim& \frac{\alpha+1}{\sqrt{2 \pi \, \alpha}} n^{-\frac{1}{2}} \left(\frac{n}{\frac{\alpha+1}{\alpha}\beta}\right)^{\beta} \left(\frac{n}{(\alpha+1)(n-\beta)}\right)^{n-\beta}
\, (1+\alpha)^{n} \notag \\ [1mm]
&\sim& \frac{\alpha+1}{\sqrt{2 \pi \, \alpha}} \, n^{-\frac{1}{2}} \, (1+\alpha)^{n} \exp\left( \beta \log\left( \frac{n}{\frac{\alpha+1}{\alpha}\beta}\right) + (n-\beta) \log\left( \frac{n}{(\alpha+1)(n-\beta)} \right)\right). \label{equation:20}
\end{eqnarray}
Before continuing, we pause to note that
\begin{equation*}
n-\frac{1}{\alpha}<\frac{\alpha+1}{\alpha} \beta \leq n+1,
\end{equation*}
and thus (recalling the definition of $\beta$)
\begin{equation*}
-1\leq n- \frac{\alpha+1}{\alpha} \beta<\frac{1}{\alpha}  \setlength{\baselineskip}{18pt}.
\end{equation*}
Simplifying the above yields
\begin{equation*}
-\alpha \leq \alpha n - (\alpha+1)\beta < 1.
\end{equation*}
Then
\begin{equation}
\alpha n- (\alpha+1) \beta = \mathcal{O}(1),
\end{equation}
or, equivalently,
\begin{equation}\label{eqstar}
n-(\alpha+1)(n-\beta)= \mathcal{O}(1).
\end{equation}
By Taylor's expansion, we have
\begin{eqnarray*}
&&\beta \log\left( \frac{n}{\frac{\alpha+1}{\alpha}\beta}\right)
= \beta \log\left(1+ \frac{n-\frac{\alpha+1}{\alpha}\beta}{\frac{\alpha+1}{\alpha}\beta}\right) \\
&=& \beta \left( \frac{n-\frac{\alpha+1}{\alpha}\beta}{\frac{\alpha+1}{\alpha}\beta} - \frac{1}{2}\left(\frac{n-\frac{\alpha+1}{\alpha}\beta}{\frac{\alpha+1}{\alpha}\beta}\right)^2 + o\left( \left(\frac{n-\frac{\alpha+1}{\alpha}\beta}{\frac{\alpha+1}{\alpha}\beta}\right)^2 \right) \right) \\
&=& \frac{\alpha n - (\alpha+1)\beta}{\alpha+1} -\frac{1}{2} \, \frac{(\alpha n - (\alpha+1)\beta)^2}{(\alpha+1)^2 \beta}
+ o\left( \frac{(\alpha n - (\alpha+1)\beta)^2}{(\alpha+1)^2 \beta} \right)  \\
&=& \frac{\alpha n - (\alpha+1)\beta}{\alpha+1} -\frac{1}{2} \, \frac{\mathcal{O}(1)}{\mathcal{O}(n)} + o\left(\frac{\mathcal{O}(1)}{\mathcal{O}(n)}\right) \quad \text{by} \,\, (\ref{eqstar}) \\
&=& \frac{\alpha n - (\alpha+1)\beta}{\alpha+1} + \mathcal{O}(n^{-1}).
\end{eqnarray*}
Similarly,
\begin{eqnarray*}
&&(n-\beta) \log\left( \frac{n}{(\alpha+1)(n-\beta)} \right)
= \frac{n-(\alpha+1)(n-\beta)}{\alpha+1} + \mathcal{O}(n^{-1}).
\end{eqnarray*}
Thus
\begin{eqnarray*}
&&\beta \log\left( \frac{n}{\frac{\alpha+1}{\alpha}\beta}\right) +
(n-\beta) \log\left( \frac{n}{(\alpha+1)(n-\beta)} \right) \\
&=& \frac{\alpha n - (\alpha+1)\beta}{\alpha+1} + \mathcal{O}(n^{-1}) + \frac{n-(\alpha+1)(n-\beta)}{\alpha+1} + \mathcal{O}(n^{-1}) \\
&=& \mathcal{O}(n^{-1}).
\end{eqnarray*}
The Proposition now follows by plugging in the above result into \eqref{equation:20}.
\end{proof}

\begin{proposition} \label{proposition6}
For sufficiently large $n$, we have
\begin{equation}
 {n \choose \big[ \frac{\alpha(n+1)}{\alpha +1}- \sqrt{n} \big]}\alpha^{\big[ \frac{\alpha(n+1)}{\alpha +1}-\sqrt{n} \big]}
 = \left( \frac{\alpha+1}{\sqrt{2 \pi \, \alpha}} \exp\left(-\frac{(\alpha+1)^2}{2\alpha}\right) + o(1) \right) n^{-\frac{1}{2}} (1+\alpha)^{n}.
\end{equation}
\end{proposition}

\begin{proof}
Let $\gamma$ denote $\big[ \frac{\alpha(n+1)}{\alpha +1} - \sqrt{n}\,\big] $, then it follows easily that, for sufficiently large $n$,
\begin{eqnarray*}
 \gamma &=& \left( \frac{\alpha}{\alpha +1} +o(1) \right) n,  \\
  n-\gamma &=& \left( \frac{1}{\alpha +1} +o(1) \right) n.
\end{eqnarray*}
By Stirling's formula, we have
\begin{eqnarray}
&&{n \choose \big[ \frac{\alpha(n+1)}{\alpha +1} -\sqrt{n}\big]}\alpha^{\big[ \frac{\alpha(n+1)}{\alpha +1}-\sqrt{n}\big]}
= {n \choose \gamma} \alpha^{\gamma} = \frac{n!}{\gamma! (n-\gamma)!} \alpha^{\gamma} \notag \\ [2mm]
&\sim& \frac{\sqrt{2 \pi n} \Big(\frac{n}{e}\Big)^n}{\sqrt{2 \pi \gamma} \Big(\frac{\gamma}{e}\Big)^\gamma  \sqrt{2 \pi (n-\gamma)} \Big(\frac{n-\gamma}{e}\Big)^{n-\gamma}} \, \alpha^{\gamma} \notag  \\ [1mm]
&\sim& \frac{1}{\sqrt{2\pi}} \sqrt{\frac{n}{\gamma(n-\gamma)}}\left(\frac{n}{\gamma}\right)^{\gamma} \left(\frac{n}{n-\gamma}\right)^{n-\gamma} \, \alpha^{\gamma} \notag \\ [1mm]
&\sim& \frac{\alpha+1}{\sqrt{2 \pi \, \alpha}} n^{-\frac{1}{2}} \left(\frac{n}{\gamma}\right)^{\gamma} \left(\frac{n}{n-\gamma}\right)^{n-\gamma} \, \alpha^{\gamma} \notag \\ [1mm]
&\sim& \frac{\alpha+1}{\sqrt{2 \pi \, \alpha}} n^{-\frac{1}{2}} \left(\frac{n}{\frac{\alpha+1}{\alpha}\gamma}\right)^{\gamma} \left(\frac{n}{(\alpha+1)(n-\gamma)}\right)^{n-\gamma}
\, (1+\alpha)^{n} \notag \\ [1mm]
&\sim& \frac{\alpha+1}{\sqrt{2 \pi \, \alpha}} \, n^{-\frac{1}{2}} \, (1+\alpha)^{n} \exp\left( \gamma \log\left( \frac{n}{\frac{\alpha+1}{\alpha}\gamma}\right) + (n-\gamma) \log\left( \frac{n}{(\alpha+1)(n-\gamma)} \right)\right). \label{equation:21}
\end{eqnarray}
By the definition of $\gamma$, we have
\begin{equation*}
\frac{\alpha(n+1)}{\alpha+1} - \sqrt{n}-1 < \gamma \leq \frac{\alpha(n+1)}{\alpha+1} - \sqrt{n}.
\end{equation*}
Hence,
\begin{equation*}
 (\alpha +1)\sqrt{n} -\alpha \leq \alpha n - (\alpha+1)\gamma < (\alpha +1)\sqrt{n}+1,
\end{equation*}
which implies that
\begin{equation*}
\alpha n - (\alpha+1)\gamma = (\alpha+1 +o(1))\sqrt{n}.
\end{equation*}
To assess (\ref{equation:21}) we first note by Taylor's expansion that
\begin{eqnarray*}
&&\gamma \log\left( \frac{n}{\frac{\alpha+1}{\alpha}\gamma}\right)
= \gamma \log\left(1+ \frac{n-\frac{\alpha+1}{\alpha}\gamma}{\frac{\alpha+1}{\alpha}\gamma}\right) \\ [1mm]
&=& \gamma \left( \frac{n-\frac{\alpha+1}{\alpha}\gamma}{\frac{\alpha+1}{\alpha}\gamma} - \frac{1}{2}\left(\frac{n-\frac{\alpha+1}{\alpha}\gamma}{\frac{\alpha+1}{\alpha}\gamma}\right)^2 + o\left( \left(\frac{n-\frac{\alpha+1}{\alpha}\gamma}{\frac{\alpha+1}{\alpha}\gamma}\right)^2 \right) \right) \\ [1mm]
&=& \frac{\alpha n - (\alpha+1)\gamma}{\alpha+1} -\frac{1}{2} \, \frac{(\alpha n - (\alpha+1)\gamma)^2}{(\alpha+1)^2 \gamma}
+ o\left( \frac{(\alpha n - (\alpha+1)\gamma)^2}{(\alpha+1)^2 \gamma} \right)  \\ [1mm]
&=& \frac{\alpha n - (\alpha+1)\gamma}{\alpha+1} -\frac{1}{2} \, \frac{\big((\alpha+1 + o(1))\sqrt{n}\big)^2}{(\alpha+1)^2 \big(\frac{\alpha}{\alpha+1}+o(1)\big)n} + o\left(\frac{\big((\alpha+1 + o(1))\sqrt{n}\big)^2}{(\alpha+1)^2 \big(\frac{\alpha}{\alpha+1}+o(1)\big)n}\right) \\
&=& \frac{\alpha n - (\alpha+1)\gamma}{\alpha+1} -\frac{1}{2}\, \frac{\alpha+1}{\alpha} + o(1).
\end{eqnarray*}
Similarly,
\begin{eqnarray*}
&&(n-\gamma) \log\left( \frac{n}{(\alpha+1)(n-\gamma)} \right)
= \frac{n-(\alpha+1)(n-\gamma)}{\alpha+1} - \frac{1}{2} (\alpha +1) + o(1).
\end{eqnarray*}
Thus
\begin{eqnarray*}
&&\gamma \log\left( \frac{n}{\frac{\alpha+1}{\alpha}\gamma}\right) + (n-\gamma) \log\left( \frac{n}{(\alpha+1)(n-\gamma)} \right) \\ [1mm]
&=& \frac{\alpha n - (\alpha+1)\gamma}{\alpha+1} - \frac{1}{2} \, \frac{\alpha+1}{\alpha} + o(1)
+ \frac{n-(\alpha+1)(n-\gamma)}{\alpha+1} - \frac{1}{2} (\alpha +1) + o(1) \\[1mm]
&=& -\frac{1}{2} \, \frac{(\alpha+1)^2}{\alpha} +o(1).
\end{eqnarray*}
Plugging in the above result into \eqref{equation:21}, Proposition \ref{proposition6} follows.
\end{proof}

With the above propositions in hand, we now turn towards the finiteness of
$\varphi_d(1-)$.

\begin{proposition} \label{proposition7}
Let $d$ be integer satisfying $d \geq 3$. For sufficiently large $n$, we have
\begin{eqnarray}
 \sum_{k=0}^{n} {n \choose k}^{d} P^{k} &\leq& \left(\frac{P^{\frac{1}{d}}+1}{\sqrt{2\pi\, P^{\frac{1}{d}}}} +o(1) \right)^{d-1} n^{-\frac{d-1}{2}} \left( 1+P^{\frac{1}{d}} \right)^{dn}, \label{equation:22} \\ [1mm]
 \sum_{k=0}^{n} {n \choose k}^{d} P^{k} &\geq& \left(\frac{P^{\frac{1}{d}}+1}{\sqrt{2\pi\, P^{\frac{1}{d}}}} \exp\left( - \frac{\left( P^{\frac{1}{d}} +1 \right)^2}{2 P^{\frac{1}{d}}}\right) +o(1) \right)^{d} n^{-\frac{d-1}{2}} \left( 1+P^{\frac{1}{d}} \right)^{dn}. \label{equation:23}
\end{eqnarray}
\end{proposition}

\begin{proof}
Set
\begin{equation*}
 \beta = \left[\frac{P^{\frac{1}{d}}(n+1)}{P^{\frac{1}{d}}+1}\right], \
 \text{and} \ \gamma = \left[\frac{P^{\frac{1}{d}}(n+1)}{P^{\frac{1}{d}}+1}-\sqrt{n}\right].
\end{equation*}
By Proposition \ref{proposition4}, we have
\begin{eqnarray*}
 {n \choose k} \left(P^{\frac{1}{d}}\right)^k &\leq& {n \choose \beta} \left(P^{\frac{1}{d}}\right)^{\beta}
 , \quad k \in \{0,1,\cdots,n\},  \\ [1mm]
 {n \choose k} \left(P^{\frac{1}{d}}\right)^k &\geq& {n \choose \gamma} \left(P^{\frac{1}{d}}\right)^{\gamma}
 , \quad k \in \{\gamma, \gamma+1, \cdots, \beta\}.
\end{eqnarray*}
Hence,
\begin{eqnarray}
\sum_{k=0}^{n} {n \choose k}^{d} P^{k}
&=& \sum_{k=0}^{n} \left({n \choose k} \left(P^{\frac{1}{d}}\right)^{k} \right)^{d-1} \,{n \choose k} \left(P^{\frac{1}{d}}\right)^{k} \notag \\[1mm]
&\leq& \sum_{k=0}^{n} \left({n \choose \beta} \left(P^{\frac{1}{d}}\right)^{\beta} \right)^{d-1} \,{n \choose k} \left(P^{\frac{1}{d}}\right)^{k}  \notag \\ [1mm]
&=& \left({n \choose \beta} \left(P^{\frac{1}{d}}\right)^{\beta} \right)^{d-1} \sum_{k=0}^{n} {n \choose k} \left(P^{\frac{1}{d}}\right)^{k} \notag \\ [1mm]
&=& \left({n \choose \beta} \left(P^{\frac{1}{d}}\right)^{\beta} \right)^{d-1} \left(1+P^{\frac{1}{d}}\right)^{n}.
 \label{equation:24}
\end{eqnarray}
By Proposition \ref{proposition5}, we have
\begin{equation*}
{n \choose \beta} \left(P^{\frac{1}{d}}\right)^{\beta} =\left(\frac{P^{\frac{1}{d}}+1}{\sqrt{2\pi\, P^{\frac{1}{d}}}} +o(1) \right)\, n^{-\frac{1}{2}} \left( 1+ P^{\frac{1}{d}}\right)^{n}.
\end{equation*}
The inequality in \eqref{equation:22} now follows by plugging the above result into \eqref{equation:24}.  Further,
\begin{eqnarray}
\sum_{k=0}^{n} {n \choose k}^{d} P^{k}
&=& \sum_{k=0}^{n} \left({n \choose k} \left(P^{\frac{1}{d}} \right)^{k} \right)^{d} \notag \\ [1mm]
&\geq& \sum_{k=\gamma}^{\beta} \left({n \choose k} \left(P^{\frac{1}{d}} \right)^{k} \right)^{d} \notag \\ [1mm]
&\geq& \sum_{k=\gamma}^{\beta} \left({n \choose \gamma} \left(P^{\frac{1}{d}} \right)^{\gamma} \right)^{d}
\notag \\ [1mm]
&=& (\beta-\gamma+1)  \left({n \choose \gamma} \left(P^{\frac{1}{d}} \right)^{\gamma} \right)^{d}. \label{equation:25}
\end{eqnarray}
By Proposition \ref{proposition6}, we have
\begin{equation*}
{n \choose \gamma} \left(P^{\frac{1}{d}} \right)^{\gamma} =
\left(\frac{P^{\frac{1}{d}}+1}{\sqrt{2\pi\, P^{\frac{1}{d}}}} \exp\left( - \frac{\left( P^{\frac{1}{d}} +1 \right)^2}{2 P^{\frac{1}{d}}}\right) +o(1) \right)\,n^{-\frac{1}{2}} \left( 1+ P^{\frac{1}{d}}\right)^{n}.
\end{equation*}
Together with \eqref{equation:25} and the fact that
\begin{eqnarray*}
&&\beta-\gamma+1= \left[\frac{P^{\frac{1}{d}}(n+1)}{P^{\frac{1}{d}}+1}\right]
- \left[\frac{P^{\frac{1}{d}}(n+1)}{P^{\frac{1}{d}}+1}-\sqrt{n}\right] +1 \\ [1mm]
&>& \frac{P^{\frac{1}{d}}(n+1)}{P^{\frac{1}{d}}+1} -1 - \left(\frac{P^{\frac{1}{d}}(n+1)}{P^{\frac{1}{d}}+1}-\sqrt{n}\right) +1 \\ [1mm]
&=&\sqrt{n},
\end{eqnarray*}
the inequality in \eqref{equation:23} follows.
\end{proof}

Proposition \ref{proposition6} tells us that $Q^n \sum_{k=0}^{n} {n \choose k}^d P^k $ has same order as $n^{- (d-1)/2} \left( Q\left(1+P^{1/d} \right)^d \right)^n$. Our next goal is to determine the value of $  Q\left(1+P^{1/d} \right)^d $. By the definition of $P$ and $Q$, we have
\begin{equation}
Q\left(1+P^{\frac{1}{d}} \right)^d
= \prod_{j=1}^{d} q_{j} \left( 1 + \left(\prod_{j=1}^{d} p_j q_j^{-1}\right)^{\frac{1}{d}} \right)^d
= \left( \left(\prod_{j=1}^{d} p_j\right)^{\frac{1}{d}} + \left(\prod_{j=1}^{d} q_j\right)^{\frac{1}{d}} \right)^d.
\label{equation:26}
\end{equation}

\begin{proposition} \label{proposition8}
\begin{equation}
\left(\prod_{j=1}^{d} p_j\right)^{\frac{1}{d}} + \left(\prod_{j=1}^{d} q_j\right)^{\frac{1}{d}} \leq 1, \label{equation:27}
\end{equation}
where equality holds if and only if $p_1 = \cdots = p_d$.
\end{proposition}

\begin{proof}
Since $f(x)=\log x$ is concave, we have
\begin{eqnarray}
 \log \left(\prod_{j=1}^{d} q_j\right)^{\frac{1}{d}} = \frac{1}{d} \sum_{j=1}^{d} \log q_j
 \leq \log \left( \frac{1}{d} \sum_{j=1}^{d} q_j \right)
 = \log  \left( 1- \frac{1}{d} \sum_{j=1}^d p_j \right). \label{equation:28}
\end{eqnarray}
Note that the inequality of arithmetic and geometric means implies $\frac{1}{d} \sum_{j=1}^d p_j \geq \left(\prod_{j=1}^{d} p_j\right)^{\frac{1}{d}}$. Together with \eqref{equation:28}, we have
\begin{equation*}
\log \left(\prod_{j=1}^{d} q_j\right)^{\frac{1}{d}}
\leq \log \left( 1- \left(\prod_{j=1}^{d} q_j\right)^{\frac{1}{d}} \right),
\end{equation*}
which implies \eqref{equation:27}. The equality holds only if $\frac{1}{d} \sum_{j=1}^d p_j = \left(\prod_{j=1}^{d} p_j\right)^{\frac{1}{d}}$, namely, $p_1 =\cdots = p_d$. If $p_1 =\cdots = p_d$ the equality holds trivially. This completes the proof.
\end{proof}
Combining Proposition \ref{proposition7}, \ref{proposition8} and equations \eqref{equation:16}and \eqref{equation:26}, Theorem \ref{theorem:3} below follows immediately.

\begin{theorem} \label{theorem:3}
In the case $d=3$, i.e. three i.i.d. random walks which are independent of each other, if $p_1=p_2=p_3$, then $\varphi_{3}(1-)=\infty$, which means $P\left(J^3 = \infty\right) =0$, i.e. rencontre happens almost surely; if $p_1,p_2,p_3$ are not equal, then $\varphi_{3}(1-)<\infty$, which means $P\left(J^3 = \infty\right) >0$. In the case $d\geq 4$, $\varphi_{d}(1-)<\infty$ regardless of the values of $p_1,\cdots,p_d$. This means that $P\left(J^d = \infty\right) >0$.
\end{theorem}

As promised in the introduction, we now provide an alternative proof that the expectation of $J^d$ is infinite.
\begin{theorem} \label{theorem:4}
For $d\geq 3$,  $E\left(J^d\right) = \infty$.
\end{theorem}
\begin{proof}
According to Theorem \ref{theorem:3}, we only need prove $E\left(J^d\right) = \infty$ in the case that $d=3$ and $p_1=p_2=p_3$, since in other cases, $P\left(J^d = \infty\right) >0$, which implies immediately that $E\left(J^d\right)=\infty$. If so, $P\left(J^3 = \infty\right) =0$, and hence
\begin{equation}
E\left(J^3\right) = \sum_{n=1}^{\infty} n \, P\left(J^3=n\right). \label{equation:30}
\end{equation}
Note that $\phi_{3}(x)$ and $\varphi_{3}(x)$ are analytic if $x \in [0,1)$. By Abel's theorem for power series, we have
\begin{eqnarray*}
&&\sum_{n=1}^{\infty} n \, P\left(J^3=n\right)
= \lim_{x \rightarrow 1-} \sum_{n=1}^{\infty} n \, P\left(J^3=n\right) x^{n-1} \\
&=& \lim_{x \rightarrow 1-} \phi_3'(x)
= \lim_{x \rightarrow 1-} \left( 1- \frac{1}{1+ \varphi_3(x)} \right)'
= \lim_{x \rightarrow 1-} \frac{\varphi'_3(x)}{\left( 1+ \varphi_3(x) \right)^2}.
\end{eqnarray*}
Together with \eqref{equation:30}, we obtain
\begin{eqnarray}
E\left(J^3\right) = \lim_{x \rightarrow 1-} \frac{\varphi'_3(x)}{\left( 1+ \varphi_3(x) \right)^2}. \label{equation:31}
\end{eqnarray}
We thus need only estimate $ {\varphi'_3(x)}/{\left( 1+ \varphi_3(x) \right)^2} $. To do so, we need introduce further notation. Let
\begin{eqnarray*}
 K_1 & :=& \left(\frac{P^{\frac{1}{3}}+1}{\sqrt{2\pi\, P^{\frac{1}{3}}}}  \right)^{2}  \\
 K_2 & :=& \left(\frac{P^{\frac{1}{3}}+1}{\sqrt{2\pi\, P^{\frac{1}{3}}}} \exp\left( - \frac{\left( P^{\frac{1}{3}} +1 \right)^2}{2 P^{\frac{1}{3}}}\right)  \right)^{3}  \\
 T & :=& Q\left(1+P^{\frac{1}{3}}\right)^{3} = \left( \left(p_1\,p_2\,p_3\right)^{\frac{1}{3}} +\left(q_1\,q_2\,q_3\right)^{\frac{1}{3}} \right)^{3} .
\end{eqnarray*}
By Proposition \ref{proposition8}, in the case  $d=3$ and $p_1=p_2=p_3$, we have $T=1$. Now consider  Proposition \ref{proposition7} with $d=3$. There exists an integer $N$ such that for $n \geq N$,
\begin{eqnarray*}
\sum_{k=0}^{n} {n \choose k}^{3} P^{k} &\leq& 2K_1 \,n^{-1} \left( 1+P^{\frac{1}{3}} \right)^{3n}, \\
\sum_{k=0}^{n} {n \choose k}^{3} P^{k} &\geq& \frac{K_2}{2} \,n^{-1} \left( 1+P^{\frac{1}{3}} \right)^{3n}.
\end{eqnarray*}
From \eqref{eqforphi}, for $0 \leq x < 1$, we have
\begin{eqnarray}
\varphi_{3}'(x) &=& \sum_{n=1}^{\infty} n\,Q^n x^{n-1} \sum_{k=0}^{n} {n \choose k}^{3} P^{k}
\geq  \sum_{n=N}^{\infty} n\,Q^n x^{n-1} \sum_{k=0}^{n} {n \choose k}^{3} P^{k} \notag \\
&\geq& \sum_{n=N}^{\infty} n\,Q^n x^{n-1} \cdot \frac{K_2}{2} \,n^{-1} \left( 1+P^{\frac{1}{3}} \right)^{3n}
= \frac{K_2}{2} \sum_{n=N}^{\infty}  \frac{1}{x} \left(Tx\right)^{n} \notag  \\
&=& \frac{K_2}{2} \sum_{n=N}^{\infty} \frac{1}{x} \, x^{n}
= \frac{K_2}{2} \frac{x^{N-1}}{1-x}. \label{equation:32}
\end{eqnarray}
Recalling Taylor's expansion for $-\log(1-x) = \sum_{n=1}^{\infty} x^{n}/n$, for $0 \leq x < 1$,
\begin{eqnarray}
1+\varphi_3(x) &=& 1+ \sum_{n=1}^{\infty} Q^n x^{n} \sum_{k=0}^{n} {n \choose k}^{3} P^{k} \notag  \\
&=& 1+ \sum_{n=1}^{N-1} Q^n x^{n} \sum_{k=0}^{n} {n \choose k}^{3} P^{k}
+ \sum_{n=N}^{\infty} Q^n x^{n} \sum_{k=0}^{n} {n \choose k}^{3} P^{k} \notag \\
&\leq& 1+ \sum_{n=1}^{N-1} Q^{n} \sum_{k=0}^{n} {n \choose k}^{3} P^{k}
+ \sum_{n=N}^{\infty} Q^n x^{n} \cdot 2K_1 \,n^{-1} \left( 1+P^{\frac{1}{3}} \right)^{3n} \notag \\
&=&  1+ \sum_{n=1}^{N-1} Q^{n} \sum_{k=0}^{n} {n \choose k}^{3} P^{k}
+ 2K_1 \sum_{n=N}^{\infty} \frac{1}{n} \left(Tx\right)^{n}  \notag   \\
&=&  1+ \sum_{n=1}^{N-1} Q^{n} \sum_{k=0}^{n} {n \choose k}^{3} P^{k}
+ 2K_1 \sum_{n=N}^{\infty} \frac{1}{n} x^{n}  \notag \\
&\leq& 1+ \sum_{n=1}^{N-1} Q^{n} \sum_{k=0}^{n} {n \choose k}^{3} P^{k}
+ 2K_1 \sum_{n=1}^{\infty} \frac{1}{n} x^{n} \notag  \\
&=& 1+ \sum_{n=1}^{N-1} Q^{n} \sum_{k=0}^{n} {n \choose k}^{3} P^{k}
- 2K_1 \log(1-x).  \label{equation:33}
\end{eqnarray}
Let
\begin{equation*}
K_3 := 1 + \sum_{n=1}^{N-1} Q^n \sum_{k=0}^{n} {n \choose k }^3 P^{k}.
\end{equation*}
From \eqref{equation:33},
\begin{equation}
1+\varphi_3(x) \leq K_3 - 2K_1 \log(1-x). \label{equation:34}
\end{equation}
Combining \eqref{equation:31}, \eqref{equation:32} and \eqref{equation:34} yields
\begin{eqnarray*}
 &&E\left( J^3\right) = \lim_{x \rightarrow 1-} \frac{\varphi_{3}'(x)}{(1+ \varphi_{3}(x))^{2}}
 \geq \lim_{x \rightarrow 1-} \frac{ \frac{K_2}{2} \frac{x^{N-1}}{1-x}}{\left( K_3- 2K_1 \log(1-x)\right)^2}
 = \infty,
\end{eqnarray*}
completing the proof.
\end{proof}

\begin{remark}
The apt referee has pointed out that the dependence of dimension $d$ in Theorems \ref{theorem:3} and \ref{theorem:4} is somewhat reminiscent of that of P\'olya's theorem for simple random walks (see \cite{Doyle}).
\end{remark}

\section{Conditional expected first rencontre-time}\label{sec6}

As we have seen throughout the preceding sections, rencontres are typically rare events. In fact, we know that $E(J^d)=\infty$ for $d \ge 2$, and $P(J^d=\infty)>0$ for all $d > 3$. Still, even rare events do happen, and of course there are many examples in science where it was the occurrence of a rare event that has given rise to new questions. However, in many of these examples, the questions are difficult to answer, in particular since they are of the a-posteriori type. A well-known example of such a question is as follows: we are here, and thus life exists, but then how plausible is it that life was born at random out of chaos?\\
\indent One way to approach such questions is to consider a system is determined by $c$ components, of which $c-1$ are assumed known and the remaining one is unknown. One may then attempt plausibility arguments for the last component to have functioned in one way or another such that the event which we see could have occurred. Our focus here is related to such objectives, although on a much more modest level. \\
\indent Specifically, suppose that $d=3, p_1=.3,$ and $p_2=.5$, and that $p:=p_3$ is unknown. The larger $p$ becomes, the more likely it is that $S^{3}(n)$ will quickly dominate $S^{1}(n)$ and $S^{2}(n)$, and so by the law of large numbers, a rencontre after time $n$ tends quickly to zero as $n$ becomes large. In other words, by knowing $J^d<\infty$ and $E(J^d | J^d<\infty)=t$, we would expect $p$ to be larger as $t$ becomes smaller because the conditional probabilities of $J^d$ given $J^d<\infty$ must be more concentrated on the smaller values of $J^d$. Our approach will be simpler in the sense that we will not work with partially unknown parameters; we instead suppose that all parameters are known and develop tools to provide bounds for 
$E(J^d | J^d<\infty)$. With $p_1$ and $p_2$ fixed, we obtain a ``sampled'' version of what we want by plugging in several values of $p_3$.\\r

With this motivation in hand, we now consider the problem raised in the introduction of calculating the conditional expectations $E(J^d|J^d<\infty)$ and $E(J^d|b < J^d<\infty)$. To obtain the bounds needed for these conditional expectations, we shall replace Stirling's formula by Robbins version of Stirling's formula:
for $n \in \mathds{N}_{+}$,
\begin{equation}
\sqrt{2 \pi}\, n^{n + \frac{1}{2}}\, e^{-n} \,e^{\frac{1}{12n+1}}
\leq n!
\leq \sqrt{2 \pi}\, n^{n + \frac{1}{2}}\, e^{-n} \,e^{\frac{1}{12n}}  \label{eq:stirling}.
\end{equation}

We shall first extend Proposition \ref{proposition5} and \ref{proposition6}. It is assumed throughout that $\alpha$ is positive. As above, the notation $[x]$ is used to denote the largest integer which is less than or equal to $x$.

\begin{proposition}  \label{prop:upperboundforalpha}
Let $\lambda$ be a real number in $(0,1)$ and let $N(\alpha, \lambda) = \max \{ [\alpha/\lambda] +1, [1/(\lambda\alpha)] +1  \}$. For $n \geq N(\alpha, \lambda)$, we have
\begin{equation*}
{n \choose \big[ \frac{\alpha(n+1)}{\alpha+1} \big] } \alpha^{\big[ \frac{\alpha(n+1)}{\alpha+1} \big] }
\leq \frac{M(\alpha, \lambda)}{\sqrt{2 \pi}} \, n^{-\frac{1}{2}}\, (1+\alpha)^{n},
\end{equation*}
where
\begin{equation*}
 M(\alpha, \lambda) := \frac{\alpha +1}{\sqrt{\alpha}} \, \frac{1}{1-\lambda}.
\end{equation*}
\end{proposition}

\begin{proof}
See Appendix B.
\end{proof}

\begin{proposition} \label{prop:lowerboundyamma}
Let $\lambda$ be a real number in $(0,1)$. If $n$ satisfies $\lambda \alpha n - (\alpha +1) \sqrt{n} -1 \geq 0$, then
\begin{eqnarray*}
&&{n \choose \big[ \frac{\alpha(n+1)}{\alpha+1} -\sqrt{n} \big] } \alpha^{\big[ \frac{\alpha(n+1)}{\alpha+1} -\sqrt{n} \big] }
\geq \frac{\sqrt{2 \pi}}{e^2} \, C_{1}(\alpha, \lambda) \, n^{-\frac{1}{2}}\, (1+\alpha)^{n},
\end{eqnarray*}
where
\begin{eqnarray*}
C_{1}(\alpha, \lambda) := \max\left\{ 4,  \frac{(\alpha+1)^2}{\alpha(1+\lambda \alpha)} \right\} \cdot
\exp\left( - \frac{1}{2} \,\frac{1+\lambda\alpha}{(1-\lambda)\alpha}\left( (\alpha+1)^2 + \frac{\lambda\alpha}{1+ \lambda\alpha} \right) \right).
\end{eqnarray*}
\end{proposition}

\begin{proof}
See Appendix B.
\end{proof}

\begin{proposition} \label{prop:lowerboundyammatilde}
Let $\lambda$ be a real number in $(0,1)$. If positive integer $n$ satisfies $\lambda  n - (\alpha +1) \sqrt{n} -\alpha \geq 0$, then
\begin{eqnarray*}
&&{n \choose \big[ \frac{\alpha(n+1)}{\alpha+1} + \sqrt{n} \big] } \alpha^{\big[ \frac{\alpha(n+1)}{\alpha+1} +\sqrt{n} \big] }
\geq \frac{\sqrt{2 \pi}}{e^2} \, C_{2}(\alpha, \lambda) \, n^{-\frac{1}{2}}\, (1+\alpha)^{n},
\end{eqnarray*}
where
\begin{eqnarray*}
C_{2}(\alpha, \lambda) := \max\left\{ 4,  \frac{(\alpha+1)^2}{\alpha+\lambda} \right\} \cdot
\exp\left( - \frac{1}{2} \,\frac{\alpha + \lambda}{(1-\lambda)\alpha^2}\left( (\alpha+1)^2 + \frac{\lambda\alpha^2}{ \alpha + \lambda} \right) \right).
\end{eqnarray*}
\end{proposition}

\begin{proof}
See Appendix B.
\end{proof}
With the above propositions in hand, we now give bounds for the coefficients of $\varphi_{d}(x)$.
\begin{proposition} \label{prop:upperboundforcoefficient}
Let $d \geq 3$ be an integer. For $n\geq N\left(P^{\frac{1}{d}}, \lambda\right),$ with $N(\alpha, \lambda)$ defined in Proposition \ref{prop:upperboundforalpha},
\begin{equation}
\sum_{k=0}^{n} {n \choose k}^{d} \, P^k
\leq \left(\frac{M\left(P^{\frac{1}{d}} , \lambda\right)}{\sqrt{2 \, \pi}}\right)^{d-1} \, n^{-\frac{d-1}{2}}
\left(1 + P^{\frac{1}{d}}\right)^{dn}. \label{eq:upperboundforcoefficients}
\end{equation}
\end{proposition}
\begin{proof}
See Appendix B.
\end{proof}

\begin{proposition}  \label{prop:lowerboundforcoefficient}
Let $d \geq 3$ be a positive integer. Define
\begin{eqnarray*}
L(\alpha, \lambda)
&:=&\max \Bigg\{  \left[ \left( \frac{(\alpha +1) + \sqrt{(\alpha +1)^2 + 4 \lambda \alpha}}{2 \lambda \alpha} \right)^2 \right] +1 ,\\ & & \qquad \quad \left[ \left( \frac{(\alpha +1) + \sqrt{(\alpha +1)^2 + 4 \lambda \alpha}}{2 \lambda } \right)^2 \right] +1 \Bigg\}.
\end{eqnarray*}
For $n \geq L\left( P^{\frac{1}{d}}, \lambda \right)$, we have
\begin{eqnarray}
\sum_{k=0}^{n} {n \choose k}^{d} P^{k} \geq \left( \frac{\sqrt{2\,\pi}}{e^2} \right)^{d}\,K\left(P^{\frac{1}{d}}, d, \lambda \right)\,n^{-\frac{d-1}{2}}\left( 1+P^{\frac{1}{d}} \right)^{dn},
\label{eq:lowerboundforcoefficient}
\end{eqnarray}
with $K(\alpha, d, \lambda)$ defined as
\begin{equation*}
K(\alpha, d, \lambda):= \frac{(1-\lambda)\alpha+1}{\alpha +1} \,\left( C_{1}(\alpha, \lambda) \right)^{d} +\left( C_{2}(\alpha, \lambda) \right)^{d},
\end{equation*}
and $C_{1}(\alpha, \lambda)$ and $C_{2}(\alpha, \lambda)$ defined, respectively, in Proposition \ref{prop:lowerboundyamma} and Proposition \ref{prop:lowerboundyammatilde}.
\end{proposition}

\begin{proof}
See Appendix B.
\end{proof}

\subsection{Bounds for the generating function}
With the above propositions in hand, we now give bounds for $\varphi_{d} (x)$ and $\varphi'_{d}(x)$. It follows from \eqref{equation:7} and \eqref{eqforphi} that
\begin{equation*}
\varphi_{d}(x) = \sum_{n=1}^{\infty} P\left(R_{n}^{d}\right) x^{n}
= \sum_{n=1}^{\infty} x^n Q^{n} \sum_{k=0}^{n} {n \choose k}^{d} P^{k},
\end{equation*}
with
\begin{equation*}
P\left(R_{n}^{d}\right) = Q^{n} \sum_{k=0}^{n} {n \choose k}^{d} P^{k}, \quad \text{for $n \in \mathds{N}_{+}$}.
\end{equation*}
Applying Proposition \ref{prop:upperboundforcoefficient} yields
\begin{eqnarray}
&&\varphi_{d}(x) = \sum_{n=1}^{N(P^{\frac{1}{d}}, \lambda )-1} P\left(R_{n}^{d}\right) x^{n}
     + \sum_{n = N(P^{\frac{1}{d}}, \lambda )}^{\infty} x^n Q^{n} \sum_{k=0}^{n} {n \choose k}^{d} P^{k} \notag  \\  [1mm]
&\leq& \sum_{n=1}^{N(P^{\frac{1}{d}}, \lambda )-1}  P\left(R_{n}^{d}\right) x^n
 \, + \sum_{n = N(P^{\frac{1}{d}}, \lambda )}^{\infty} x^n Q^{n} \left( \frac{M\left(P^{\frac{1}{d}},\lambda\right)}{\sqrt{2\,\pi}} \right)^{d-1} \, n^{- \frac{d-1}{2}} \,
 \left(1 + P^{\frac{1}{d}}\right)^{dn}   \notag  \\  [1mm]
&=&  \sum_{n=1}^{N(P^{\frac{1}{d}}, \lambda )-1} P\left(R_{n}^{d}\right) x^n  + \left( \frac{M\left(P^{\frac{1}{d}},\lambda\right)}{\sqrt{2\,\pi}} \right)^{d-1} \sum_{n = N(P^{\frac{1}{d}}, \lambda )}^{\infty} n^{- \frac{d-1}{2}} \, \left( Q\left(1 + P^{\frac{1}{d}}\right)^{d}\, x \right)^{n}.  \label{eq:upperforvarphi}
\end{eqnarray}
Let $UB\left(x;P,Q,d, \lambda | \varphi_{d}\right)$ denote the right-hand side of \eqref{eq:upperforvarphi}, i.e. the upper bound for $\varphi_{d} (x)$.
Applying Proposition \ref{prop:lowerboundforcoefficient} to $\varphi_{d}(x)$ yields
\begin{eqnarray}
&&\varphi_{d}(x) = \sum_{n=1}^{L(P^{\frac{1}{d}}, \lambda )-1} P\left(R_{n}^{d}\right) x^{n} +
 \sum_{n = L(P^{\frac{1}{d}}, \lambda )}^{\infty}  x^n  Q^{n} \sum_{k=0}^{n} {n \choose k}^{d} P^{k}   \notag \\
&\geq&  \sum_{n=1}^{L(P^{\frac{1}{d}}, \lambda )-1} P\left(R_{n}^{d}\right) x^{n} +
  \sum_{n = L(P^{\frac{1}{d}}, \lambda )}^{\infty}  x^n  Q^{n} \left(\frac{\sqrt{2\,\pi}}{e^2}\right)^{d}
K\left( P^{\frac{1}{d}}, d,  \lambda \right) n^{-\frac{d-1}{2}} \left(1+P^{\frac{1}{d}}\right)^{dn} \notag \\
&=& \sum_{n=1}^{L(P^{\frac{1}{d}}, \lambda )-1} P\left(R_{n}^{d}\right) x^{n} + \left(\frac{\sqrt{2\,\pi}}{e^2}\right)^{d}
K\left( P^{\frac{1}{d}},d, \lambda \right)   \sum_{n = L(P^{\frac{1}{d}}, \lambda )}^{\infty}
n^{- \frac{d-1}{2}} \, \left( Q\left(1 + P^{\frac{1}{d}}\right)^{d}\, x \right)^{n}.  \label{eq:lowerforvarphi}
\end{eqnarray}
Let $LB\left(x;P,Q,d, \lambda | \varphi_{d}\right)$ denote the right-hand side of \eqref{eq:lowerforvarphi}, i.e. the lower bound for $\varphi_{d} (x)$. It follows easily from \eqref{eq:upperforvarphi} that $\varphi_{d}(x)$ is convergent for $0 \leq x < \Big( Q\big( 1 + P^{\frac{1}{d}} \big)^{d} \Big)^{-1}$, and hence $\varphi_{d}(x)$ is analytic in this region. Then
\begin{equation*}
\varphi'_{d} (x) = \sum_{n=1}^{\infty} n \, P\left(R_{n}^{d}\right) x^{n-1}
= \sum_{n=1}^{\infty} n \, x^{n-1} Q^{n} \sum_{k=0}^{n} {n \choose k}^{d} P^{k},
\end{equation*}
and
\begin{equation*}
\varphi''_{d} (x) = \sum_{n=1}^{\infty} n(n-1) \, P\left(R_{n}^{d}\right) x^{n-2}
= \sum_{n=1}^{\infty} n(n-1) \, x^{n-2} Q^{n} \sum_{k=0}^{n} {n \choose k}^{d} P^{k}.
\end{equation*}
Similarly, applying Proposition \ref{prop:upperboundforcoefficient} and \ref{prop:lowerboundforcoefficient} to $\varphi'_{d}(x)$, we have for $0 < x < \Big( Q\big( 1 + P^{\frac{1}{d}} \big)^{d} \Big)^{-1}$,
\begin{equation}
LB\left(x;P,Q,d, \lambda | \varphi'_{d}\right)
\leq \varphi'_{d}(x)
\leq UB\left(x;P,Q,d, \lambda | \varphi'_{d}\right),  \label{eq:boundforvarphiprime}
\end{equation}
where $UB\left(x;P,Q,d, \lambda | \varphi'_{d}\right)$ is defined as
\begin{eqnarray}
\sum_{n=1}^{N(P^{\frac{1}{d}}, \lambda )-1} n\, P\left(R^{d}_{n}\right) \, x^{n-1}  + \left( \frac{M\left(P^{\frac{1}{d}},\lambda\right)}{\sqrt{2\,\pi}} \right)^{d-1} \sum_{n = N(P^{\frac{1}{d}}, \lambda )}^{\infty} n^{- \frac{d-3}{2}}\,x^{-1} \, \left( Q\left(1 + P^{\frac{1}{d}}\right)^{d}\, x \right)^{n}, \label{eq:upperforvarphiprime}
\end{eqnarray}
and $LB\left(x;P,Q,d, \lambda | \varphi'_{d}\right)$ is defined as
\begin{eqnarray}
\sum_{n=1}^{L(P^{\frac{1}{d}}, \lambda )-1} n\, P\left(R^{d}_{n}\right) \, x^{n-1} +
\left(\frac{\sqrt{2\,\pi}}{e^2}\right)^{d}
K\left( P^{\frac{1}{d}}, d, \lambda \right)   \sum_{n = L(P^{\frac{1}{d}}, \lambda )}^{\infty}
n^{- \frac{d-3}{2}}\, x^{-1} \, \left( Q\left(1 + P^{\frac{1}{d}}\right)^{d}\, x \right)^{n}.  \label{eq:lowerforvarphiprime}
\end{eqnarray}
Note that
\begin{equation*}
\varphi'_{d}(x) + x \varphi''_{d}(x)  = \sum_{n=1}^{\infty} n^2 \, x^{n-1} Q^{n} \sum_{k=0}^{n} {n \choose k}^{d} P^{k}.
\end{equation*}
Applying Proposition \ref{prop:upperboundforcoefficient} to $\varphi'_{d}(x) + x \varphi''_{d}(x)$, we have for $0 < x < \Big( Q\big( 1 + P^{\frac{1}{d}} \big)^{d} \Big)^{-1}$,
\begin{equation}
\varphi'_{d}(x) + x \varphi''_{d}(x)
\leq UB\left(x;P,Q,d, \lambda | \varphi'_{d} + x \varphi''_{d}\right),  \label{eq:boundforvarphisecondprime}
\end{equation}
where $UB\left(x;P,Q,d, \lambda | \varphi'_{d} + x \varphi'' \right)$ is defined as
\begin{eqnarray}
\sum_{n=1}^{N(P^{\frac{1}{d}}, \lambda )-1} n^2\, P\left(R_{n}^{d}\right) \, x^{n-1}  + \left( \frac{M\left(P^{\frac{1}{d}},\lambda\right)}{\sqrt{2\,\pi}} \right)^{d-1} \sum_{n = N(P^{\frac{1}{d}}, \lambda )}^{\infty} n^{- \frac{d-5}{2}}\,x^{-1} \, \left( Q\left(1 + P^{\frac{1}{d}}\right)^{d}\, x \right)^{n}. \label{eq:upperforvarphisecondprime}
\end{eqnarray}

\subsection{Bounds for $E\left( J^{d} | J^{d} < \infty \right)$}
Recall that in Section \ref{sec5}, we have shown the expected value of $J^{d}$ to always be infinite (see Theorem \ref{theorem:4}). We now investigate the conditional expectation $E\left( J^{d} | J^{d} < \infty \right)$ and give bounds for it.

We first observe that
\begin{equation*}
E\left( J^{d} | J^{d} < \infty \right) = \frac{\sum_{n=1}^{\infty} n\,P\left(J^d = n\right) }{\sum_{n=1}^{\infty} P\left(J^d = n\right)} = \frac{\lim_{x \rightarrow 1-} \phi'_{d}(x)}{ \lim_{x \rightarrow 1-} \phi_{d}(x)}
= \lim_{x \rightarrow 1-} \frac{\phi'_{d}(x)}{\phi_{d} (x)}.
\end{equation*}
The last equality holds because the limit of $\phi_{d}(x)$ is positive and finite as $x$ tends to $1-$. Since $\phi_{d}(x) = 1- \frac{1}{1 + \varphi_{d}(x)}$ (i.e. Lemma \ref{proposition:1}), we have
\begin{equation}
E\left( J^{d} | J^{d} < \infty \right) = \lim_{x \rightarrow 1-} \frac{\phi'_{d}(x)}{\phi_{d} (x)}
= \lim_{x \rightarrow 1-} \frac{\varphi'_{d}(x)}{\varphi_{d}(x) \left(1+ \varphi_{d}(x)\right)}. \label{eq:representationforfirstmoment}
\end{equation}
Applying the bounds for $\varphi_{d}(x)$ and $\varphi'_{d}(x)$, i.e. \eqref{eq:upperforvarphi}, \eqref{eq:lowerforvarphi}, \eqref{eq:boundforvarphiprime}, \eqref{eq:upperforvarphiprime} and \eqref{eq:lowerforvarphiprime}, with $\lambda$ replaced by $\lambda_1$ in the upper bounds and $\lambda$  replaced by $\lambda_2$ in the lower bounds, Theorem \ref{thm:boundforconditionalexpectationinfty} immediately follows.
\begin{theorem}  \label{thm:boundforconditionalexpectationinfty}
Let $d \geq 3$ be a positive integer and let $\lambda_1$ and $\lambda_2$ be two arbitrary real numbers in $(0,1)$. We have
\begin{eqnarray}
E\left( J^{d} | J^{d} < \infty \right) &\leq& \lim_{x \rightarrow 1-} \frac{UB\left(x;P,Q,d, \lambda_1 | \varphi'_{d}\right)}{LB\left(x;P,Q,d, \lambda_2 | \varphi_{d}\right) \left(1+LB\left(x;P,Q,d, \lambda_2 | \varphi_{d}\right)\right)},  \label{eq:upperboundforconditionalexpectationinfty}  \\  [1mm]
E\left( J^{d} | J^{d} < \infty \right) &\geq& \lim_{x \rightarrow 1-} \frac{LB\left(x;P,Q,d, \lambda_2 | \varphi'_{d}\right)}{UB\left(x;P,Q,d, \lambda_1 | \varphi_{d}\right) \left(1+UB\left(x;P,Q,d, \lambda_1 | \varphi_{d}\right)\right)}. \label{eq:lowerboundforconditionalexpectationinfty}
\end{eqnarray}
\end{theorem}

If $Q\big( 1 + P^{\frac{1}{d}} \big)^{d} < 1$ (i.e. $p_1, \cdots , p_d$ are not all the same, see Proposition \ref{proposition8} and \eqref{equation:26}), then $UB\left(1;P,Q,d, \lambda_1 | \varphi'_{d}\right)$ and $LB\left(1;P,Q,d, \lambda_2 | \varphi_{d}\right)$ are both finite, since the power series in \eqref{eq:lowerforvarphi} and \eqref{eq:upperforvarphiprime} are convergent when $x=1$. Hence, by \eqref{eq:upperboundforconditionalexpectationinfty}, $E\left( J^{d} | J^{d} < \infty \right)$ is also finite. Note that if $Q\big( 1 + P^{\frac{1}{d}} \big)^{d} = 1$ (i.e. $p_1=\cdots = p_d$) and $d = 4$ or $5$, then the power series in \eqref{eq:lowerforvarphiprime} diverges when $x=1$ 
but the power series in \eqref{eq:upperforvarphi} converges when $x=1$,
i.e. $LB\left(1;P,Q,d, \lambda_2 | \varphi'_{d}\right) = \infty$ but $UB\left(1;P,Q,d, \lambda_1 | \varphi_{d}\right) < \infty$. In this case, it follows immediately from \eqref{eq:lowerboundforconditionalexpectationinfty} that $E\left( J^{d} | J^{d} < \infty \right) = \infty$. If $Q\big( 1 + P^{\frac{1}{d}} \big)^{d} = 1$ and $d \geq 6$, note that the series $\sum_{n} n^{- a}$ converges for $a >1$, and thus $UB\left(1;P,Q,d, \lambda_1 | \varphi'_{d}\right)$ and $LB\left(1;P,Q,d, \lambda_2 | \varphi_{d}\right)$ are both finite. Hence $E\left( J^{d} | J^{d} < \infty \right)$ is  again finite by \eqref{eq:upperboundforconditionalexpectationinfty}. The only remaining case to consider is $Q\big( 1 + P^{\frac{1}{d}} \big)^{d} = 1$ and $d=3$. In this case, as $x \rightarrow 1-$,
\begin{eqnarray*}
&&LB\left(x;P,Q,d, \lambda_2 | \varphi'_{d}\right)
\geq \left(\frac{\sqrt{2\,\pi}}{e^2}\right)^{3} K\left( P^{\frac{1}{3}},3, \lambda_2 \right)
\sum_{n= L\big(P^{\frac{1}{3}}, \lambda_2\big)}^{\infty} x^{-1} \cdot x^{n}  \\
&=& \left(\frac{\sqrt{2\,\pi}}{e^2}\right)^{3} K\left( P^{\frac{1}{3}},3, \lambda_2 \right)
  \cdot \frac{x^{L(P^{\frac{1}{3}}, \lambda_2)-1}}{1-x}
= \mathcal{O}\left((1-x)^{-1}\right),
\end{eqnarray*}
and
\begin{eqnarray*}
UB\left(x;P,Q,d, \lambda_1 | \varphi_{d}\right) &\leq&
\sum^{N(P^{\frac{1}{3}}, \lambda_1 ) -1 }_{n=1} x^{n}  + \left( \frac{M\left(P^{\frac{1}{3}},\lambda_1\right)}{\sqrt{2\,\pi}} \right)^{2} \sum_{n = N(P^{\frac{1}{3}}, \lambda_1 )}^{\infty} n^{- 1} \, x^{n}  \\
&\leq& \sum^{N(P^{\frac{1}{3}}, \lambda_1 ) -1 }_{n=1} 1  + \left( \frac{M\left(P^{\frac{1}{3}},\lambda_1\right)}{\sqrt{2\,\pi}} \right)^{2} \sum_{n = 1}^{\infty} n^{- 1} \, x^{n}  \\
&=&  N\left(P^{\frac{1}{3}}, \lambda_1\right)-1  - \left( \frac{M\left(P^{\frac{1}{3}},\lambda_1\right)}{\sqrt{2\,\pi}} \right)^{2} \log(1-x) \\
&=& \mathcal{O}\left(-\log(1-x)\right),
\end{eqnarray*}
The order of the numerator of the right-hand side of \eqref{eq:lowerboundforconditionalexpectationinfty} is at least $\mathcal{O}\left((1-x)^{-1}\right)$ but the order of the denominator is at most $\mathcal{O}\left( \left(\log(1-x)\right)^{2}  \right)$ as $x \rightarrow 1-$, which implies the right-hand side of \eqref{eq:lowerboundforconditionalexpectationinfty} tends to $\infty$ as $x$ tends to $1-$. Hence, $E\left( J^{d} | J^{d} < \infty \right) = \infty$.

The above results are now summarized by Corollary \ref{coro:conditionalisinfty} below.
\begin{corollary} \label{coro:conditionalisinfty}
Let $d\geq 3$ be a positive integer, then
\begin{equation*}
E\left( J^{d} | J^{d} < \infty \right)
\begin{cases}
= \infty, & \text{if $p_1 = \cdots = p_d$ and $d \in \{3,4,5\}$,}  \\
< \infty, & \text{otherwise.}
\end{cases}
\end{equation*}
\end{corollary}
We conclude by offering numerics of the bounds for $E\left(J^d | J^d < \infty\right)$ in Table \ref{tab:numericalresult11} below.

\begin{table}[!htbp]
\centering

\begin{tabular}{|l|c|c|}
\hline
Parameter settings & Lower bound  & Upper bound   \\
\hline
\tabincell{l}{$d=3, p_1=.3, p_2=.4, p_3=.5$ \\ $\lambda_1=1/80, \lambda_2 = 1/8$}  & 3.86223 & 3.88172 \\
\hline
\tabincell{l}{$d=3, p_1=.45, p_2=.5, p_3=.55$ \\ $\lambda_1=1/300, \lambda_2 = 1/10$}  & 9.31034 & 9.84928 \\
\hline
\tabincell{l}{$d=3, p_1=.05, p_2=.5, p_3=.5$ \\ $\lambda_1=1/15, \lambda_2 = 1/2$}  & 1.22586 & 1.22586 \\
\hline
\tabincell{l}{$d=4, p_1=.3, p_2=.4, p_3=.5, p_4=.6$ \\ $\lambda_1=1/15, \lambda_2 = 1/2$}  & 2.3814 & 2.38296 \\
\hline
\tabincell{l}{$d=4, p_1=.4, p_2=.45, p_3=.5, p_4=.55$ \\ $\lambda_1=1/250, \lambda_2 = 1/8$}  & 4.35938 & 4.361 \\
\hline
\tabincell{l}{$d=4, p_1=.47, p_2=.5, p_3=.52, p_4=.53$ \\ $\lambda_1=1/500, \lambda_2 = 1/15$}  & 9.9011 & 10.3937 \\
\hline
\tabincell{l}{$d=4, p_1=.5, p_2=.5, p_3=.6, p_4=.6$ \\ $\lambda_1=1/200, \lambda_2 = 1/8$}  & 4.73906 & 4.75067 \\
\hline
\tabincell{l}{$d=4, p_1=.48, p_2=.49, p_3=.5, p_4=.51$ \\ $p_5 = 0.52$, $\lambda_1=1/500, \lambda_2 = 1/15$}  & 5.1569 & 5.49917 \\
\hline
\tabincell{l}{$d=4, p_1=.4, p_2=.4, p_3=.5, p_4=.5$ \\ $p_5 = 0.5$, $\lambda_1=1/150, \lambda_2 = 1/8$}  & 3.02342 & 3.0273 \\
\hline
\end{tabular}
\caption{Numerics for upper and lower bounds of $E(J^d | J^{d} < \infty)$} \label{tab:numericalresult11}.
\end{table}

\subsection{Bounds for $E(J^d | b < J^d < \infty )$}

We shall now find a upper bound for $E\left(J^d | b < J^d < \infty \right)$ for small $b$. For ease of notation, let us define a new random variable $\tilde{J}^{d}$ to be a positive-integer-valued random variable equaling $n$ with probability $P\left(J^d =n\right)/P\left(J^d < \infty\right)$. That is, $\tilde{J}^{d}$ is $J^d$ conditioned on the event $\left\{J^d < \infty\right\}$. As such, $E\left(J^d| J^d < \infty\right) = E\big( \tilde{J}^{d} \big)$. We shall henceforth let $\mu$ denote the expectation of $\tilde{J}^{d}$.

\begin{theorem}
Let $d \geq 3$ be a positive integer and let $t$ be a positive real number in $(1,\infty)$. Let $\lambda_1$ and $\lambda_2$ be arbitrary real numbers in $(0,1)$. If $p_1, \cdots , p_d$ are not all the same or $d \geq 6$, then
\begin{eqnarray}
&& E\left(J^d \Big| \frac{\mu}{t} < J^d < \infty \right) \notag \\
 &\leq&  \frac{t^2}{(t-1)^{2}}
\cdot \lim_{x \rightarrow 1-} \left( \frac{UB\left(x;P,Q,d, \lambda_1 | \varphi'_{d} + x \varphi''_{d}\right)}{LB\left(x;P,Q,d, \lambda_2 | \varphi'_{d}\right)} - \frac{LB\left(x;P,Q,d, \lambda_2 | \varphi'_{d}\right)}{1+ UB\left(x;P,Q,d, \lambda_1 | \varphi_{d}\right)} \right).
\end{eqnarray}
\end{theorem}

\begin{proof}
By the definition of conditional expectation and the definition of $\tilde{J}^{d}$, we have
\begin{eqnarray}
&&E\left(J^d \Big| \frac{\mu}{t} < J^d < \infty \right)
= \frac{\sum_{n= [\mu/t]+1}^{\infty} n \, P\left(J^{d}=n\right) }{\sum_{n= [\mu/t]+1}^{\infty}  P\left(J^{d}=n\right)} \notag \\ [1mm]
&=& \frac{\sum_{n= [\mu/t]+1}^{\infty} n \, P\left(J^{d}=n\right) }{\sum_{n= 1}^{\infty}  P\left(J^{d}=n\right)} \Bigg/
\frac{\sum_{n= [\mu/t]+1}^{\infty}  P\left(J^{d}=n\right) }{\sum_{n= 1}^{\infty}  P\left(J^{d}=n\right)} \notag  \\ [1mm]
&\leq& \frac{\sum_{n= 1}^{\infty} n \, P\left(J^{d}=n\right) }{\sum_{n= 1}^{\infty}  P\left(J^{d}=n\right)} \Bigg/
\frac{\sum_{n= [\mu/t]+1}^{\infty}  P\left(J^{d}=n\right) }{\sum_{n= 1}^{\infty}  P\left(J^{d}=n\right)} \notag   \\ [1mm]
&=& E\left( \tilde{J}^{d} \right) \Big/ P\left( \tilde{J}^{d} > \mu /t \right).  \label{eq:6.1}
\end{eqnarray}
By the conditional form of Jensen's inequality,
\begin{eqnarray*}
&& E\left[ \left(\tilde{J}^{d}\right)^{2} \right] \geq P\left( \tilde{J}^{d} > \mu /t \right) \, E\left[ \left(\tilde{J}^{d}\right)^{2} \bigg| \tilde{J}^{d} > \mu/t \right]   \\ [1mm]
&\geq&  P\left( \tilde{J}^{d} > \mu /t \right) \, \left[E\left( \tilde{J}^{d} \Big| \tilde{J}^{d} > \mu/t \right) \right]^2   \\  [1mm]
&=& \frac{ \left[ E\left( \tilde{J}^d ; \tilde{J}^{d} > \mu/t \right) \right]^2 }{P\left( \tilde{J}^{d} > \mu /t \right)}
= \frac{ \left[ E\left(\tilde{J}^{d}\right) - E\left( \tilde{J}^d ; \tilde{J}^{d} \leq \mu/t \right) \right]^2 }{P\left( \tilde{J}^{d} > \mu /t \right)}  \\ [1mm]
&\geq& \frac{\left( \mu - \mu/t \right)^2}{ P\left( \tilde{J}^{d} > \mu /t \right) }
= \frac{(t-1)^2}{t^2} \, \frac{\mu^2}{P\left( \tilde{J}^{d} > \mu /t \right)}.
\end{eqnarray*}
Thus,
\begin{equation*}
P\left( \tilde{J}^{d} > \mu /t \right) \geq \frac{(t-1)^2}{t^2} \, \frac{\mu^2}{E\left[ \left(\tilde{J}^{d}\right)^{2} \right]}.
\end{equation*}
Together with \eqref{eq:6.1} and the fact that $E\big(\tilde{J}^{d}\big)=\mu$, it follows that
\begin{eqnarray}
E\left(J^d \Big| \frac{\mu}{t} < J^d < \infty \right)
\leq \frac{t^2}{(t-1)^2 \, \mu} \,E\left[ \left(\tilde{J}^{d}\right)^{2} \right]. \label{eq:6.2}
\end{eqnarray}
We now represent the right-hand side of \eqref{eq:6.2} in terms of $\varphi_{d}(x)$ and its derivatives.
\begin{eqnarray*}
E\left[ \left(\tilde{J}^{d}\right)^{2} \right] &=& \sum_{n=1}^{\infty} n^2 \, P\left( \tilde{J}^{d} =n \right)
= \sum_{n=1}^{\infty} n(n-1) \, P\left( \tilde{J}^{d} =n \right) + \sum_{n=1}^{\infty} n \, P\left( \tilde{J}^{d} =n \right) \\  [1mm]
&=& \frac{\sum_{n=1}^{\infty} n(n-1) \, P\left( J^{d} =n \right)}{\sum_{n=1}^{\infty} P\left( J^{d} =n \right)}
 + \frac{\sum_{n=1}^{\infty} n \, P\left( J^{d} =n \right)}{\sum_{n=1}^{\infty} P\left( J^{d} =n \right)} \\ [1mm]
&=& \frac{\lim_{x \rightarrow 1-} \phi''_{d}(x) }{\lim_{x \rightarrow 1-} \phi_{d}(x) }
  + \frac{\lim_{x \rightarrow 1-} \phi'_{d}(x) }{\lim_{x \rightarrow 1-} \phi_{d}(x) }
  = \lim_{x \rightarrow 1-} \frac{\phi'_{d}(x) + \phi''_{d}(x)}{\phi_{d}(x)},
\end{eqnarray*}
where the last step follows since $\lim_{x \rightarrow 1-} \phi_{d}(x)$ is finite and in $(0,1]$. Together with the fact $\phi_{d}(x) = \varphi_{d}(x)/\left( 1+ \varphi_{d}(x) \right)$, it follows that
\begin{equation}
E\left[ \left(\tilde{J}^{d}\right)^{2} \right] =
\lim_{x \rightarrow 1-} \left( \frac{\varphi''_{d}(x)}{\varphi_{d}(x) \left( 1+ \varphi_{d}(x) \right) }
+  \frac{\varphi'_{d}(x)}{\varphi_{d}(x) \left( 1+ \varphi_{d}(x) \right) }
- \frac{2 \cdot \left(\varphi'_{d}(x)\right)^2}{\varphi_{d}(x) \left( 1+ \varphi_{d}(x) \right)^2 }
\right). \label{eq:representationforsecondmoment}
\end{equation}
We know from \eqref{eq:representationforfirstmoment}) that
\begin{equation*}
\mu = E\left( J^d | J^{d} < \infty \right) = \lim_{x \rightarrow 1-} \frac{\varphi'_{d}(x)}{\varphi_{d}(x) \left( 1+ \varphi_{d}(x) \right) }.
\end{equation*}
Since $\mu$ is finite and positive, we can interchange the orders of the limit and the fraction. Hence,
\begin{eqnarray*}
&&E\left[ \left(\tilde{J}^{d}\right)^{2} \right] \bigg/ \mu
= \lim_{x \rightarrow 1-} \left( \frac{\varphi''_{d}(x)}{\varphi'_{d}(x) } + 1 -\frac{2\varphi'_{d}(x)}{ 1+ \varphi_{d}(x)  } \right) \\
&=& \lim_{x \rightarrow 1-} \left( x \cdot\frac{\varphi''_{d}(x)}{\varphi'_{d}(x) } + 1 -\frac{2\varphi'_{d}(x)}{ 1+ \varphi_{d}(x)  } \right)
= \lim_{x \rightarrow 1-} \left( \frac{\varphi'_{d}(x) + x\varphi''_{d}(x)}{\varphi'_{d}(x) } -\frac{2\varphi'_{d}(x)}{ 1+ \varphi_{d}(x)  } \right).
\end{eqnarray*}
Combining the above result with \eqref{eq:6.2} yields
\begin{equation}
E\left(J^d \Big| \frac{\mu}{t} < J^d < \infty \right)
\leq \frac{t^2}{(t-1)^2} \cdot \lim_{x \rightarrow 1-} \left( \frac{\varphi'_{d}(x) + x\varphi''_{d}(x)}{\varphi'_{d}(x) } -\frac{2\varphi'_{d}(x)}{ 1+ \varphi_{d}(x)  } \right). \label{eq:CEfor>b upperbound}
\end{equation}
Applying bounds \eqref{eq:upperforvarphi}, \eqref{eq:boundforvarphiprime}, and \eqref{eq:boundforvarphisecondprime} to \eqref{eq:CEfor>b upperbound}, and replacing $\lambda$ by $\lambda_1$ in the upper bounds and replacing $\lambda$ by $\lambda_2$ in lower bounds, the proof is completed.
\end{proof}

\section{Acknowledgments}
We thank the two anonymous referees for their invaluable comments, which have greatly improved the quality of this paper. We are grateful to Professor Larry Shepp for bringing this problem to our attention, who in turn learned of this problem from Professor Abram Kagan.

\bibliographystyle{apalike}
\bibliography{Reference}

\newpage

\section{Appendix A}

We derive an explicit expression for $1+ \varphi_{2}(x)$ (see \eqref{expforvarphi2}). For $d =2$, it follows from \eqref{equation:13} that
\begin{equation*}
    1 + \varphi_2(x)
    = \frac{1}{(2\pi)^{2}} \int_{[-\pi,\pi]^2} \frac{1}{(1-x \psi_2(\theta_1, \theta_2))(1-e^{-i (\theta_1 +\theta_2)} \, x \psi_{2}(\theta_1, \theta_2))}\  d\theta_1 d\theta_2,
\end{equation*}
where $\psi_{2}(\theta_1, \theta_2)=\left(p_1 e^{i \theta_1} + q_1\right)\left(p_2 e^{i\theta_2} + q_2\right)$. Let $z_1 = e^{i \theta_1}$ and $z_2 = e^{i \theta_2}$. 
This yields
\begin{equation*}
 1 + \varphi_2(x) = \frac{1}{(2 \pi i)^2} \oint_{\gamma \times \gamma} \frac{dz_1 dz_2}{\left( 1 - x(p_1 z_1 +q _1) (p_2 z_2 +q _2) \right)\left( z_1 z_2 - x (p_1 z_1 + q_1) (p_2 z_2 +q _2) \right)},
\end{equation*}
where $\gamma$ is a counter-clockwise unit circle with center at $0$. We first calculate the integral with respect to $z_1$.
Let
\begin{eqnarray*}
 A_1 &=& 1- x q_1 (p_2 z_2 + q_2),  \\
 B_1 &=& x p_1 (p_2 z_2 + q_2),   \\
 C_1 &=& z_2 - x p_1 (p_2 z_2 + q_2), \\
 D_1 &=& x q_1 (p_2 z_2 + q_2).
\end{eqnarray*}
Then
\begin{eqnarray}
 1 + \varphi_2 (x) &=& \frac{1}{(2 \pi i)^2} \oint_{\gamma^2} \frac{dz_1 dz_2}{ (A_1 - B_1z_1)(C_1 z_1 - D_1) } \notag \\
  &=& \frac{1}{(2 \pi i)^2} \oint_{\gamma^2} \frac{dz_1 dz_2}{B_1 C_1 (\frac{A_1}{B_1} - z_1)(z_1 - \frac{D_1}{C_1}) }.
  \label{equd2:1}
\end{eqnarray}
Note that since $|p_2 z_2 + q_2| \leq p_2 |z_2| + q_2 = p_2 + q_2 =1 $,
\begin{eqnarray*}
 &&|A_1| = |1- x q_1 (p_2 z_2 + q_2)| \geq 1- |xq_1 (p_2 z_2 + q_2)| \geq 1- x q_1 > x-x q_1 = x p_1 \geq |B_1|, \\
 &&|C_1| = |z_2 - x p_1 (p_2 z_2 + q_2)| \geq |z_2| - |x p_1 (p_2 z_2 + q_2)| \geq 1 - xp_1>x - xp_1 = xq_1 \geq |D_1|.
\end{eqnarray*}
This implies that $\left| \frac{A_1}{B_1} \right|>1$ and $\left| \frac{D_1}{C_1} \right|<1$. Hence, $\frac{1}{\frac{A_1}{B_1} - z_1}$ is analytic in the unit disk and $ \frac{1}{z_1 - \frac{D_1}{C_1}}$ has a simple pole at $z_1 = \frac{D_1}{C_1}$ in unit disk. The integral \eqref{equd2:1} may then be calculated as
\begin{eqnarray}
 &&1+ \varphi_2(x) = \frac{1}{2 \pi i} \oint_{\gamma} \frac{dz_2}{ B_1 C_1 (\frac{A_1}{B_1}- \frac{D_1}{C_1}) }
 = \frac{1}{2 \pi i} \oint_{\gamma} \frac{dz_2}{ A_1 C_1 - B_1 D_1 } \notag  \\
 &=& \frac{1}{2 \pi i} \oint_{\gamma} \frac{dz_2} {\left(1- x q_1 (p_2 z_2 + q_2)\right) \left(z_2 - x p_1 (p_2 z_2 + q_2)\right) - x p_1 (p_2 z_2 + q_2) x q_1 (p_2 z_2 + q_2) } \notag \\
 &=& \frac{1}{2 \pi i} \oint_{\gamma} \frac{dz_2}{-xq_1p_2z_2^2 + (1-xp_1p_2-xq_1q_2)z_2 - x p_1q_2} \label{equd2:2}.
\end{eqnarray}
Now let
\begin{eqnarray*}
w_1 &=& \frac{(1-xp_1p_2 - xq_1q_2) + \sqrt{1-2x(p_1p_2 + q_1q_2) + x^2(p_1p_2 - q_1q_2)^2}}{2xq_1p_2}, \\
w_2 &=& \frac{(1-xp_1p_2 - xq_1q_2) - \sqrt{1-2x(p_1p_2 + q_1q_2) + x^2(p_1p_2 - q_1q_2)^2}}{2xq_1p_2},
\end{eqnarray*}
i.e., $w_1$ and $w_2$ are two roots of equation $-xq_1p_2z_2^2 + (1-xp_1p_2-xq_1q_2)z_2 - x p_1q_2=0$. In order for $w_1$ and $w_2$ to be well defined, we need to show that
\begin{equation}
1-2x(p_1p_2 + q_1q_2) + x^2(p_1p_2 - q_1q_2)^2 \geq (1-x)^2 >0, \qquad \text{for $0<x<1$}.  \label{equd2:3}
\end{equation}
Let $s_j = p_j-q_j$, $j=1,2$, then $|s_j| \leq 1$, and $p_j = (1+s_j)/2$, $q_j=(1-s_j)/2$. Then
\begin{eqnarray*}
p_1 p_2 - q_1 q_2 &=& \frac{1+s_1}{2}\cdot\frac{1+s_2}{2} - \frac{1-s_1}{2}\cdot\frac{1-s_2}{2} =\frac{s_1+s_2}{2}, \\
p_1 p_2 + q_1 q_2 &=& \frac{1+s_1}{2}\cdot\frac{1+s_2}{2} + \frac{1-s_1}{2}\cdot\frac{1-s_2}{2} =\frac{1+s_1s_2}{2}.
\end{eqnarray*}
This implies that
\begin{eqnarray}
 && 1-2x(p_1p_2 + q_1q_2) + x^2(p_1p_2 - q_1q_2)^2 \notag \\
 &=& 1- 2x\,\frac{1+s_1s_2}{2} + x^2 \, \left(\frac{s_1+s_2}{2}\right)^2 \notag  \\
 &=& 1- x(1+s_1 s_2 ) + s_1 s_2 x^{2} + x^2 \, \left(\frac{s_1+s_2}{2}\right)^2 - s_1 s_2 x^{2} \notag \\
 &=& (1-x)(1-s_1 s_2 x) + \left( \frac{s_1-s_2}{2} \right)^{2} x^{2}. \label{equd2:4}
\end{eqnarray}
Then \eqref{equd2:3} follows since $x<1$ and $|s_j| \leq1$. With roots $w_1$ and $w_2$, \eqref{equd2:2} can be represented as
\begin{eqnarray}
1 + \varphi_2(x) = \frac{1}{2 \pi i} \oint_{\gamma} \frac{dz_2}{-2xq_1 p _2 (z_2 - w_1)(z_2 - w_2)}. \label{equd2:5}
\end{eqnarray}
We proceed to calculate
\begin{eqnarray*}
&& 1- xp_1 p_2 - xq_1 q_2  \\
&=& 1-x + x - xp_1 p_2 - xq_1 q_2 \\
&=& 1-x + x(p_1 +q_1)(p_2+q_2) - xp_1 p_2 - xq_1 q_2 \\
&=& 1-x + xp_1q_2 + xq_1 p_2 \\
&>& xp_1q_2 + xq_1 p_2 \\
&=& 2xq_1 p_2 + x(p_1q_2-q_1 p_2) \\ [1mm]
&=& 2xq_1 p_2 + x \left( \frac{1+s_1}{2} \cdot \frac{1-s_2}{2} - \frac{1-s_1}{2} \cdot \frac{1+s_2}{2}  \right) \\ [1mm]
&=& 2xq_1 p_2 + \left(\frac{s_1 -s_2}{2}\right)x.
\end{eqnarray*}
Combining the above result with \eqref{equd2:4}, we have
\begin{equation}
1- xp_1 p_2 - xq_1 q_2 + \sqrt{1-2x(p_1p_2 + q_1q_2) + x^2(p_1p_2 - q_1q_2)^2} > 2xq_1 p_2 \label{equd2:6}.
\end{equation}
Similarly, we have
\begin{equation}
1- xp_1 p_2 - xq_1 q_2 + \sqrt{1-2x(p_1p_2 + q_1q_2) + x^2(p_1p_2 - q_1q_2)^2} > 2xp_1 q_2 \label{equd2:7}.
\end{equation}
It follows directly from \eqref{equd2:6} that $w_1 >1$. Further note that
\begin{eqnarray*}\
w_2 &=& \frac{(1-xp_1p_2 - xq_1q_2) - \sqrt{1-2x(p_1p_2 + q_1q_2) + x^2(p_1p_2 - q_1q_2)^2}}{2xq_1p_2} \\ [1mm]
&=& \frac{2xp_1q_2}{(1-xp_1p_2 - xq_1q_2) + \sqrt{1-2x(p_1p_2 + q_1q_2) + x^2(p_1p_2 - q_1q_2)^2}}.
\end{eqnarray*}
Together with \eqref{equd2:7}, we have that $0<w_2<1$.
Then the integral in \eqref{equd2:5} can be calculated as
\begin{eqnarray*}
1+ \varphi_2(x) = \frac{1}{-2xq_1 p _2 (w_2 - w_1)} = \frac{1}{\sqrt{1-2x(p_1p_2 + q_1q_2) + x^2(p_1p_2 - q_1q_2)^2}}.
\end{eqnarray*}
It now easily follows that for $x \in [0,1)$,
\begin{equation*}
1-\phi_2(x) = \frac{1}{1+\varphi(x)}= \sqrt{1-2x(p_1p_2 + q_1q_2) + x^2(p_1p_2 - q_1q_2)^2}.
\end{equation*}

\section{Appendix B}
This appendix contains the proofs of Propositions \ref{prop:upperboundforalpha}-\ref{prop:lowerboundforcoefficient}.\\

\noindent \textbf{Proof of Proposition \ref{prop:upperboundforalpha}}

\begin{proof}
Let $\beta$ denote $\big[ \frac{\alpha(n+1)}{\alpha+1} \big]$. Then
\begin{equation*}
{n \choose \big[ \frac{\alpha(n+1)}{\alpha+1} \big] } \alpha^{\big[ \frac{\alpha(n+1)}{\alpha+1} \big] }
= {n \choose \beta} \alpha^{\beta}
= \frac{n!}{\beta! (n-\beta)!} \alpha^{\beta}.
\end{equation*}
As we shall see, under the assumption $n \geq \max \{ [\alpha/\lambda] +1, [1/(\lambda\alpha)] +1  \}$, we have
$ 1 \leq \big[\frac{\alpha(n+1)}{\alpha+1}\big] \leq n-1$, or equivalently, $1 \leq \beta \leq n-1$ and $1 \leq n-\beta \leq n-1$. Thus, applying \eqref{eq:stirling} to the above equation gives
\begin{eqnarray}
&&{n \choose \big[ \frac{\alpha(n+1)}{\alpha+1} \big] } \alpha^{\big[ \frac{\alpha(n+1)}{\alpha+1} \big] }
= \frac{n!}{\beta! (n-\beta)!} \alpha^{\beta}  \notag \\
&\leq& \frac{\sqrt{2 \pi}\, n^{n + \frac{1}{2}}\, e^{-n} \,e^{\frac{1}{12n}}}{ \left(\sqrt{2 \pi}\, \beta^{\beta + \frac{1}{2}}\, e^{-\beta} \,e^{\frac{1}{12\beta +1}}\right) \cdot \left( \sqrt{2 \pi}\, (n-\beta)^{n-\beta + \frac{1}{2}}\, e^{-(n-\beta)} \,e^{\frac{1}{12(n-\beta) +1}} \right)} \, \alpha^{\beta} \notag \\
&=& \frac{1}{\sqrt{2 \pi}} \, \left( \frac{n^2}{\beta(n-\beta)} \right)^{\frac{1}{2}}\, n^{-\frac{1}{2}}
   \, \left(\frac{n}{\frac{\alpha+1}{\alpha}\beta}\right)^{\beta} \, \left(\frac{n}{(\alpha+1)(n-\beta)}\right)^{n-\beta} \notag \\
&& \cdot \, (1+\alpha)^{n} \, e^{\frac{1}{12}\left(\frac{1}{n} - \frac{1}{\beta+\frac{1}{12}} - \frac{1}{n-\beta + \frac{1}{12}}\right)}. \label{eq:expansion}
\end{eqnarray}
Note that $f(x) = 1/x$ for $x>0$ is a convex function. By Jensen's inequality,
\begin{eqnarray*}
&& \frac{1}{\beta+\frac{1}{12}} + \frac{1}{n-\beta + \frac{1}{12}}
\geq 2 \, \frac{1}{\frac{1}{2} \left( \beta + \frac{1}{12} + n- \beta + \frac{1}{12} \right)}
 = \frac{4}{ n + \frac{1}{6}} > \frac{4}{2n} > \frac{1}{n}.
\end{eqnarray*}
This implies that
\begin{equation}
 e^{\frac{1}{12}\left(\frac{1}{n} - \frac{1}{\beta+\frac{1}{12}} - \frac{1}{n-\beta + \frac{1}{12}}\right)} \leq 1. \label{eq:1.1}
\end{equation}
Note that $N(\alpha, \lambda) > \alpha/\lambda$ and $N(\alpha, \lambda) > 1/(\lambda\alpha)$. We then have
\begin{eqnarray*}
 &&\frac{n^2}{\beta(n-\beta)}
 = \frac{n^2}{\big[\frac{\alpha(n+1)}{\alpha +1}\big]\cdot \left(n - \big[\frac{\alpha(n+1)}{\alpha +1}\big] \right)}
\leq \frac{n^2}{ \left( \frac{\alpha(n+1)}{\alpha +1} -1 \right) \cdot \left( n- \frac{\alpha(n+1)}{\alpha +1} \right) }  \\ [1mm]
&=& \frac{(1 + \alpha)^2}{\alpha} \, \frac{n^2}{(n- 1/\alpha)(n - \alpha)}
 =  \frac{(1 + \alpha)^2}{\alpha} \, \frac{n}{n-1/\alpha} \, \frac{n}{n - \alpha}  \\ [1mm]
&\leq& \frac{(1 + \alpha)^2}{\alpha} \, \frac{N(\alpha, \lambda)}{N(\alpha, \lambda)-1/\alpha} \, \frac{N(\alpha, \lambda)}{N(\alpha, \lambda) - \alpha}  \\ [1mm]
&\leq& \frac{(1 + \alpha)^2}{\alpha} \, \frac{1/(\lambda \alpha)}{1/(\lambda \alpha)-1/\alpha} \, \frac{\alpha/\lambda}{\alpha/\lambda - \alpha}
= \frac{(1 + \alpha)^2}{\alpha} \, \left( \frac{1}{1-\lambda} \right)^2.
\end{eqnarray*}
Thus,
\begin{equation}
 \left( \frac{n^2}{\beta(n-\beta)} \right)^{\frac{1}{2}} \leq M(\alpha, \lambda). \label{eq:1.2}
\end{equation}
Since the inequality $\log(1 +x) \leq x$ holds for $x > -1$, we have
\begin{equation*}
\beta \log\left(\frac{n}{\frac{\alpha+1}{\alpha} \, \beta}\right)
= \beta \log\left( 1 + \frac{n - \frac{\alpha+1}{\alpha} \, \beta}{\frac{\alpha+1}{\alpha} \, \beta}\right)
\leq \beta \cdot \frac{n - \frac{\alpha+1}{\alpha} \, \beta}{\frac{\alpha+1}{\alpha} \, \beta}
=\frac{\alpha n - (\alpha+1)\beta}{\alpha +1}.
\end{equation*}
Hence,
\begin{equation}
\left( \frac{n}{\frac{\alpha +1}{\alpha} \, \beta} \right)^{\beta}
= \exp\left\{\beta \log\left(\frac{n}{\frac{\alpha+1}{\alpha} \, \beta}\right)\right\}
\leq \exp \left\{ \frac{\alpha n - (\alpha+1)\beta}{\alpha +1} \right\}. \label{eq:1.3}
\end{equation}
Similarly,
\begin{equation}
\left( \frac{n}{(\alpha +1) (n- \beta)} \right)^{n-\beta}
\leq \exp\left\{  \frac{n- (\alpha +1)(n-\beta)}{\alpha+1} \right\}. \label{eq:1.4}
\end{equation}
Combining \eqref{eq:expansion}, \eqref{eq:1.1}, \eqref{eq:1.2}, \eqref{eq:1.3} and \eqref{eq:1.4} completes the proof.
\end{proof}

\bigskip

\noindent \textbf{Proof of Proposition \ref{prop:lowerboundyamma}}
\begin{proof}
Let $\gamma$ denote $\big[ \frac{\alpha(n+1)}{\alpha+1} -\sqrt{n} \big]$. Then
\begin{equation}  \label{nicenewequation}
{n \choose \big[ \frac{\alpha(n+1)}{\alpha+1} -\sqrt{n} \big] } \alpha^{\big[ \frac{\alpha(n+1)}{\alpha+1} -\sqrt{n} \big] }
= {n \choose \gamma} \alpha^{\gamma}
= \frac{n!}{\gamma! (n-\gamma)!} \alpha^{\gamma}.
\end{equation}
As we shall see, the assumption $\lambda \alpha n - (\alpha +1) \sqrt{n} -1 \geq 0$ ensures $1 \leq \frac{\alpha(n+1)}{\alpha +1} - \sqrt{n} < n$, hence,
$1 \leq \gamma \leq n-1$. Applying a simple bound for $n!$, i.e.
\begin{equation*}
\sqrt{2\pi} n^{n + \frac{1}{2}} e^{-n} \leq n! \leq e n^{n + \frac{1}{2}} e^{-n}
\end{equation*}
to the equation in (\ref{nicenewequation}) gives
\begin{eqnarray}
&&{n \choose \big[ \frac{\alpha(n+1)}{\alpha+1} -\sqrt{n} \big] } \alpha^{\big[ \frac{\alpha(n+1)}{\alpha+1} -\sqrt{n} \big] }
= \frac{n!}{\gamma! (n-\gamma)!} \alpha^{\gamma}  \notag \\
&\geq& \frac{\sqrt{2 \pi}\, n^{n + \frac{1}{2}}\, e^{-n} }{ \left(e\, \gamma^{\gamma + \frac{1}{2}}\, e^{-\gamma} \right) \cdot \left( e\, (n-\gamma)^{n-\gamma + \frac{1}{2}}\, e^{-(n-\gamma)}  \right)} \, \alpha^{\gamma} \notag \\
&=& \frac{\sqrt{2 \pi}}{e^2} \, \left( \frac{n^2}{\gamma(n-\gamma)} \right)^{\frac{1}{2}}\, n^{-\frac{1}{2}}
   \, \left(\frac{n}{\frac{\alpha+1}{\alpha}\gamma}\right)^{\gamma} \, \left(\frac{n}{(\alpha+1)(n-\gamma)}\right)^{n-\gamma}
 \, (1+\alpha)^{n} . \label{eq:expansion:gamma}
\end{eqnarray}
From the definition of $\gamma$, we have
\begin{equation*}
\frac{\alpha n}{\alpha +1} - \sqrt{n} - \frac{1}{\alpha +1} < \gamma
\leq \frac{\alpha n}{\alpha +1} - \sqrt{n} + \frac{\alpha}{\alpha +1},
\end{equation*}
and thus
\begin{equation*}
\frac{n}{\alpha +1} + \sqrt{n} - \frac{\alpha}{\alpha +1} \leq n- \gamma
< \frac{n}{\alpha +1} + \sqrt{n} + \frac{1}{\alpha +1}.
\end{equation*}
It follows easily from the above inequalities and by the assumption $\lambda \alpha n \geq (\alpha +1) \sqrt{n} +1 $ that
\begin{eqnarray*}
 &&\gamma (n-\gamma) < \left(\frac{\alpha n}{\alpha +1} - \sqrt{n} + \frac{\alpha}{\alpha +1}\right)
  \cdot \left( \frac{n}{\alpha +1} + \sqrt{n} + \frac{1}{\alpha +1} \right)  \\
 &\leq& \left(\frac{\alpha n}{\alpha +1}\right) \cdot \left( \frac{n}{\alpha +1} + \frac{\lambda \alpha n}{\alpha+1} \right)
 = \frac{\alpha (1+\lambda \alpha)}{(\alpha +1)^2} \, n^2.
\end{eqnarray*}
Hence
\begin{equation}
\frac{n^2}{\gamma(n-\gamma)} > \frac{(\alpha+1)^2}{\alpha(1+\lambda \alpha)} . \label{eq:2.0}
\end{equation}
By the inequality of arithmetic and geometric means, $\gamma (n - \gamma) \leq (\gamma + n -\gamma )^2/4 = n^2/4$, which implies $\frac{n^2}{\gamma(n-\gamma)} \geq 4$. Together with \eqref{eq:2.0}, we have
\begin{equation}
\frac{n^2}{\gamma(n-\gamma)} \geq \max\left\{ 4,  \frac{(\alpha+1)^2}{\alpha(1+\lambda \alpha)} \right\} . \label{eq:2.1}
\end{equation}
By the assumption $\lambda \alpha n - (\alpha +1) \sqrt{n} -1 \geq 0$, we have
\begin{equation}
\sqrt{n} \geq \frac{(\alpha+1)+\sqrt{(\alpha+1)^2 + 4 \lambda \alpha}}{2\lambda \alpha}
> \frac{\alpha+1}{\lambda \alpha}. \label{eq:boundforn}
\end{equation}
Since $\log(1+x) = \int_{0}^{x} \frac{1}{1+s} \,ds  = \int_{0}^{1} \frac{x}{1+xt} \,dt$,
\begin{eqnarray}
&& \gamma \log\left(\frac{n}{\frac{\alpha+1}{\alpha}\gamma}\right) +
 (n-\gamma) \log\left(\frac{n}{(\alpha+1)(n-\gamma)}\right) \notag  \\
&=& \gamma \log\left( 1+ \frac{\alpha n -(\alpha +1)\gamma}{(\alpha+1)\gamma} \right)
   + (n - \gamma) \log\left( 1- \frac{\alpha n -(\alpha +1)\gamma}{(\alpha+1)(n-\gamma)} \right) \notag \\ [1mm]
&=& \gamma \int_{0}^{1} \frac{\frac{\alpha n -(\alpha +1)\gamma}{(\alpha+1)\gamma}}{1 + \frac{\alpha n -(\alpha +1)\gamma}{(\alpha+1)\gamma}\,t} \, dt +
   (n-\gamma) \int_{0}^{1}  \frac{- \frac{\alpha n -(\alpha +1)\gamma}{(\alpha+1)(n-\gamma)}}{1 - \frac{\alpha n -(\alpha +1)\gamma}{(\alpha+1)(n-\gamma)} \, t} \, dt  \notag  \\
&=& - \frac{1}{(\alpha +1)^2} \frac{(\alpha n -(\alpha+1)\gamma)^2 \,n}{\gamma (n-\gamma)}
 \int_{0}^{1} \frac{t}{\left(1 + \frac{\alpha n -(\alpha +1)\gamma}{(\alpha+1)\gamma}\,t\right)  \left(1 - \frac{\alpha n -(\alpha +1)\gamma}{(\alpha+1)(n-\gamma)} \, t\right)} \, dt.  \label{eq:expressionlog}
\end{eqnarray}
One may easily verify by taking first derivatives that $\frac{(\alpha n -(\alpha+1)\gamma)^2 \,n}{\gamma (n-\gamma)}$ is a non-increasing function of $\gamma$ when $\gamma \in \big[0, \frac{\alpha n}{\alpha +1}\big]$. Hence
\begin{eqnarray}
&&\frac{(\alpha n -(\alpha+1)\gamma)^2 \,n}{\gamma (n-\gamma)} \leq
\frac{\left(\alpha n -(\alpha+1)\left(\frac{\alpha n}{\alpha +1} - \sqrt{n} - \frac{1}{\alpha +1}\right)\right)^2 \,n}{\left(\frac{\alpha n}{\alpha +1} - \sqrt{n} - \frac{1}{\alpha +1}\right) \left(n-\left(\frac{\alpha n}{\alpha +1} - \sqrt{n} - \frac{1}{\alpha +1}\right)\right)}  \notag \\  [1mm]
&=& (\alpha+1)^2  \frac{\left( (\alpha+1)\sqrt{n} +1 \right)^2 \, n}{\left(\alpha n -(\alpha +1)\sqrt{n} -1\right)\left(n+(\alpha+1)\sqrt{n} +1\right)} \notag \\ [1mm]
&\leq& (\alpha+1)^2  \frac{\left( (\alpha+1)\sqrt{n} +1 \right)^2 \, n}{(1-\lambda)\alpha n \cdot \left(n+(\alpha+1)\sqrt{n} +1\right)} \notag \\ [1mm]
&=& \frac{(\alpha+1)^2}{(1-\lambda)\alpha} \cdot \frac{\left( (\alpha+1)\sqrt{n} +1 \right)^2}{n+(\alpha+1)\sqrt{n} +1} \notag \\ [1mm]
&=& \frac{(\alpha+1)^2}{(1-\lambda)\alpha} \cdot \left(\frac{(\alpha+1)^2 n + (\alpha+1)\sqrt{n} + 1}{n+(\alpha+1)\sqrt{n} +1}  +  \frac{(\alpha+1)\sqrt{n}}{n+(\alpha+1)\sqrt{n} +1}\right) \notag \\ [1mm]
&\leq& \frac{(\alpha+1)^2}{(1-\lambda)\alpha} \cdot \left( (\alpha+1)^2 + \frac{(\alpha+1)}{\sqrt{n} + (\alpha+1)} \right)   \notag \\ [1mm]
&\leq& \frac{(\alpha+1)^2}{(1-\lambda)\alpha} \cdot \left( (\alpha+1)^2 + \frac{(\alpha+1)}{\frac{\alpha +1}{\lambda\alpha} + (\alpha+1)} \right)   \qquad \text{applying \eqref{eq:boundforn}}  \notag \\ [1mm]
&=& \frac{(\alpha+1)^2}{(1-\lambda)\alpha} \cdot \left( (\alpha+1)^2 + \frac{\lambda\alpha}{1+\lambda\alpha} \right). \label{eq:boundlog1}
\end{eqnarray}
Similarly, since $\frac{\alpha n -(\alpha+1)\gamma}{(\alpha+1)(n-\gamma)}$ is a non-increasing function of $\gamma$ when $\gamma \in (0, n)$,
\begin{eqnarray}
&&\frac{\alpha n -(\alpha+1)\gamma}{(\alpha+1)(n-\gamma)}
\leq \frac{\alpha n -(\alpha+1)\left(\frac{\alpha n}{\alpha +1} - \sqrt{n} - \frac{1}{\alpha +1}\right)}{(\alpha+1)\left(n-\left(\frac{\alpha n}{\alpha +1} - \sqrt{n} - \frac{1}{\alpha +1}\right)\right)} \notag \\
&=& \frac{(\alpha +1)\sqrt{n}+1}{n + (\alpha+1)\sqrt{n} + 1}
\leq \frac{(\alpha +1)\sqrt{n} + 1}{\frac{1}{\lambda\alpha} \, \left( (\alpha+1)\sqrt{n} + 1 \right) + (\alpha+1)\sqrt{n} + 1} = \frac{\lambda\alpha}{1+\lambda\alpha}.  \label{eq:boundlog2}
\end{eqnarray}
Combining \eqref{eq:expressionlog}, \eqref{eq:boundlog1}, \eqref{eq:boundlog2} and using the fact that $\frac{\alpha n -(\alpha+1) \gamma}{(\alpha+1)\gamma } >0$ gives us
\begin{eqnarray*}
\text{\eqref{eq:expressionlog}}
&\geq& - \frac{1}{(1-\lambda)\alpha} \left( (\alpha+1)^2 + \frac{\lambda\alpha}{1+ \lambda\alpha} \right)
 \cdot \int_{0}^{1} \frac{t}{1\cdot \left( 1- \frac{\lambda\alpha}{1+\lambda\alpha} \right)} \, dt  \\
 &=& - \frac{1}{2} \,\frac{1+\lambda\alpha}{(1-\lambda)\alpha}\left( (\alpha+1)^2 + \frac{\lambda\alpha}{1+ \lambda\alpha} \right).
\end{eqnarray*}
This implies
\begin{eqnarray}
\left(\frac{n}{\frac{\alpha+1}{\alpha}\gamma}\right)^{\gamma} \, \left(\frac{n}{(\alpha+1)(n-\gamma)}\right)^{n-\gamma}
\geq \exp\left( - \frac{1}{2} \,\frac{1+\lambda\alpha}{(1-\lambda)\alpha}\left( (\alpha+1)^2 + \frac{\lambda\alpha}{1+ \lambda\alpha} \right) \right). \label{eq:2.2}
\end{eqnarray}
Combining \eqref{eq:expansion:gamma}, \eqref{eq:2.1}, and \eqref{eq:2.2} completes the proof.
\end{proof}
\newpage
\bigskip

\noindent \textbf{Proof of Proposition \ref{prop:lowerboundyammatilde}}

\begin{proof}
Let $\tilde{\gamma}$ denote $\big[ \frac{\alpha(n+1)}{\alpha+1} + \sqrt{n} \big]$. Then
\begin{equation*}
{n \choose \big[ \frac{\alpha(n+1)}{\alpha+1} + \sqrt{n} \big] } \alpha^{\big[ \frac{\alpha(n+1)}{\alpha+1} + \sqrt{n} \big] }
= {n \choose \tilde{\gamma}} \alpha^{\tilde{\gamma}}
= \frac{n!}{\tilde{\gamma}! (n-\tilde{\gamma})!} \alpha^{\tilde{\gamma}}.
\end{equation*}
As we shall see, the assumption $\lambda  n - (\alpha +1) \sqrt{n} -\alpha \geq 0$ ensures that $1 \leq \frac{\alpha(n+1)}{\alpha +1} + \sqrt{n} < n$, and, hence,
$1 \leq \tilde{\gamma} \leq n-1$. Applying a simple bound for $n!$, i.e.
\begin{equation*}
\sqrt{2\pi} n^{n + \frac{1}{2}} e^{-n} \leq n! \leq e n^{n + \frac{1}{2}} e^{-n}
\end{equation*}
to the above equation gives
\begin{eqnarray}
&&{n \choose \big[ \frac{\alpha(n+1)}{\alpha+1} + \sqrt{n} \big] } \alpha^{\big[ \frac{\alpha(n+1)}{\alpha+1} +\sqrt{n} \big] }
= \frac{n!}{\tilde{\gamma}! (n-\tilde{\gamma})!} \alpha^{\tilde{\gamma}}  \notag \\
&\geq& \frac{\sqrt{2 \pi}\, n^{n + \frac{1}{2}}\, e^{-n} }{ \left(e\, \tilde{\gamma}^{\tilde{\gamma} + \frac{1}{2}}\, e^{-\tilde{\gamma}} \right) \cdot \left( e\, (n-\tilde{\gamma})^{n-\tilde{\gamma} + \frac{1}{2}}\, e^{-(n-\tilde{\gamma})}  \right)} \, \alpha^{\tilde{\gamma}} \notag \\
&=& \frac{\sqrt{2 \pi}}{e^2} \, \left( \frac{n^2}{\tilde{\gamma}(n-\tilde{\gamma})} \right)^{\frac{1}{2}}\, n^{-\frac{1}{2}}
   \, \left(\frac{n}{\frac{\alpha+1}{\alpha}\tilde{\gamma}}\right)^{\tilde{\gamma}} \, \left(\frac{n}{(\alpha+1)(n-\tilde{\gamma})}\right)^{n-\tilde{\gamma}}
 \, (1+\alpha)^{n} . \label{eq:expansion:gamma:tilde}
\end{eqnarray}
From the definition of $\tilde{\gamma}$, we have
\begin{equation*}
\frac{\alpha n}{\alpha +1} + \sqrt{n} - \frac{1}{\alpha +1} < \tilde{\gamma}
\leq \frac{\alpha n}{\alpha +1} + \sqrt{n} + \frac{\alpha}{\alpha +1}.
\end{equation*}
Thus,
\begin{equation*}
\frac{n}{\alpha +1} - \sqrt{n} - \frac{\alpha}{\alpha +1} \leq n- \tilde{\gamma}
< \frac{n}{\alpha +1} - \sqrt{n} + \frac{1}{\alpha +1}.
\end{equation*}
It follows easily from the above inequalities and by the assumption $\lambda n \geq (\alpha +1) \sqrt{n} +\alpha $ that
\begin{eqnarray*}
 &&\tilde{\gamma} (n-\tilde{\gamma}) < \left(\frac{\alpha n}{\alpha +1} + \sqrt{n} + \frac{\alpha}{\alpha +1}\right)
  \cdot \left( \frac{n}{\alpha +1} - \sqrt{n} + \frac{1}{\alpha +1} \right)  \\
 &\leq& \left(\frac{\alpha n}{\alpha +1} + \frac{\lambda n}{\alpha +1}\right) \cdot \left( \frac{n}{\alpha +1} \right)
 = \frac{\alpha + \lambda}{(\alpha +1)^2} \, n^2.
\end{eqnarray*}
This yields
\begin{equation}
\frac{n^2}{\tilde{\gamma}(n-\tilde{\gamma})} > \frac{(\alpha+1)^2}{\alpha+ \lambda} . \label{eq:3.0}
\end{equation}
Again, by the inequality of arithmetic and geometric means, $\tilde{\gamma} (n - \tilde{\gamma}) \leq (\tilde{\gamma} + n -\tilde{\gamma} )^2/4 = n^2/4$, which implies that $\frac{n^2}{\tilde{\gamma}(n-\tilde{\gamma})} \geq 4$. Together with \eqref{eq:3.0}, we have
\begin{equation}
\frac{n^2}{\tilde{\gamma}(n-\tilde{\gamma})} \geq \max\left\{ 4,  \frac{(\alpha+1)^2}{\alpha+ \lambda} \right\} . \label{eq:3.1}
\end{equation}
By the assumption that $\lambda  n - (\alpha +1) \sqrt{n} - \alpha \geq 0$, we have
\begin{equation}
\sqrt{n} \geq \frac{(\alpha+1)+\sqrt{(\alpha+1)^2 + 4 \lambda \alpha}}{2\lambda }
> \frac{\alpha+1}{\lambda}. \label{eq:boundforn:tilde}
\end{equation}
Since $\log(1+x) = \int_{0}^{x} \frac{1}{1+s} \,ds  = \int_{0}^{1} \frac{x}{1+xt} \,dt$,
\begin{eqnarray}
&& \tilde{\gamma} \log\left(\frac{n}{\frac{\alpha+1}{\alpha}\tilde{\gamma}}\right) +
 (n-\tilde{\gamma}) \log\left(\frac{n}{(\alpha+1)(n-\tilde{\gamma})}\right) \notag  \\
&=& \tilde{\gamma} \log\left( 1 - \frac{ (\alpha +1)\tilde{\gamma} - \alpha n}{(\alpha+1)\tilde{\gamma}} \right)
   + (n - \tilde{\gamma}) \log\left( 1+ \frac{ (\alpha +1)\tilde{\gamma} -\alpha n}{(\alpha+1)(n-\tilde{\gamma})} \right) \notag \\ [1mm]
&=& \tilde{\gamma} \int_{0}^{1} \frac{-\frac{(\alpha +1)\tilde{\gamma} -\alpha n}{(\alpha+1)\tilde{\gamma}}}{1 - \frac{(\alpha +1)\tilde{\gamma} -\alpha n}{(\alpha+1)\tilde{\gamma}}\,t} \, dt +
   (n-\tilde{\gamma}) \int_{0}^{1}  \frac{ \frac{(\alpha +1)\tilde{\gamma} -\alpha n}{(\alpha+1)(n-\tilde{\gamma})}}{1 + \frac{(\alpha +1)\tilde{\gamma} -\alpha n}{(\alpha+1)(n-\tilde{\gamma})} \, t} \, dt  \notag  \\
&=& - \frac{1}{(\alpha +1)^2} \frac{((\alpha +1)\tilde{\gamma} -\alpha n)^2 \,n}{\tilde{\gamma} (n-\tilde{\gamma})}
 \int_{0}^{1} \frac{t}{\left(1 + \frac{(\alpha +1)\tilde{\gamma} -\alpha n}{(\alpha+1)(n- \tilde{\gamma})}\,t\right)  \left(1 - \frac{(\alpha +1)\tilde{\gamma} -\alpha n}{(\alpha+1)\tilde{\gamma}} \, t\right)} \, dt.  \label{eq:expressionlog:tilde}
\end{eqnarray}
It is easy to verify by taking first derivatives that $\frac{((\alpha +1)\tilde{\gamma} -\alpha n)^2 \,n}{\tilde{\gamma} (n-\tilde{\gamma})}$ is a non-decreasing function of $\tilde{\gamma}$ when $\tilde{\gamma} \in \big[ \frac{\alpha n}{\alpha +1},n \big]$. Hence
\begin{eqnarray}
&&\frac{((\alpha +1)\tilde{\gamma} -\alpha n)^2 \,n}{\tilde{\gamma} (n-\tilde{\gamma})} \leq
\frac{\left( (\alpha+1)\left(\frac{\alpha n}{\alpha +1} + \sqrt{n} + \frac{\alpha}{\alpha +1}\right) - \alpha n\right)^2 \,n}{\left(\frac{\alpha n}{\alpha +1} + \sqrt{n} + \frac{\alpha}{\alpha +1}\right) \left(n-\left(\frac{\alpha n}{\alpha +1} + \sqrt{n} + \frac{\alpha}{\alpha +1}\right)\right)}  \notag \\  [1mm]
&=& (\alpha+1)^2  \frac{\left( (\alpha+1)\sqrt{n} + \alpha \right)^2 \, n}{\left(\alpha n + (\alpha +1)\sqrt{n} + \alpha\right)\left(n - (\alpha+1)\sqrt{n} - \alpha\right)} \notag \\ [1mm]
&\leq& (\alpha+1)^2  \frac{\left( (\alpha+1)\sqrt{n} + \alpha \right)^2 \, n}{ \left(\alpha n+(\alpha+1)\sqrt{n} + \alpha\right) \cdot (1-\lambda)n} \notag \\ [1mm]
&=& \frac{(\alpha+1)^2}{1-\lambda} \cdot \frac{\left( (\alpha+1)\sqrt{n} + \alpha \right)^2}{\alpha n+(\alpha+1)\sqrt{n} + \alpha} \notag \\ [1mm]
&=& \frac{(\alpha+1)^2}{1-\lambda} \cdot \left(\frac{(\alpha+1)^2 n + \alpha(\alpha+1)\sqrt{n} + \alpha^2}{\alpha n+(\alpha+1)\sqrt{n} + \alpha}  +  \frac{\alpha(\alpha+1)\sqrt{n}}{\alpha n+(\alpha+1)\sqrt{n} + \alpha}\right) \notag \\ [1mm]
&\leq& \frac{(\alpha+1)^2}{1-\lambda} \cdot \left( \frac{(\alpha+1)^2}{\alpha} + \frac{\alpha(\alpha+1)}{\alpha\sqrt{n} + (\alpha+1)} \right)   \notag \\ [1mm]
&\leq& \frac{(\alpha+1)^2}{1-\lambda} \cdot \left( \frac{(\alpha+1)^2}{\alpha} + \frac{\alpha(\alpha+1)}{\alpha \cdot\frac{\alpha +1}{\lambda} + (\alpha+1)} \right)   \qquad \text{apply \eqref{eq:boundforn:tilde}}  \notag \\ [1mm]
&=& \frac{(\alpha+1)^2}{(1-\lambda)\alpha} \cdot \left( (\alpha+1)^2 + \frac{\lambda\alpha^2}{\alpha + \lambda} \right). \label{eq:boundlog1:tilde}
\end{eqnarray}
Similarly, since $\frac{ (\alpha+1)\tilde{\gamma} - \alpha n}{(\alpha+1)\tilde{\gamma}}$ is a non-decreasing function of $\tilde{\gamma}$ when $\tilde{\gamma} \in (0, n)$,
\begin{eqnarray}
&&\frac{ (\alpha+1)\tilde{\gamma} - \alpha n}{(\alpha+1)\tilde{\gamma}}
\leq \frac{(\alpha+1)\left(\frac{\alpha n}{\alpha +1} + \sqrt{n} + \frac{\alpha}{\alpha +1}\right) - \alpha n}{(\alpha+1)\left(\frac{\alpha n}{\alpha +1} + \sqrt{n} + \frac{\alpha}{\alpha +1}\right)} \notag \\
&=& \frac{(\alpha +1)\sqrt{n} + \alpha}{ \alpha n + (\alpha+1)\sqrt{n} + \alpha}
\leq \frac{(\alpha +1)\sqrt{n} + \alpha}{\frac{\alpha}{\lambda} \, \left( (\alpha+1)\sqrt{n} + \alpha \right) + (\alpha+1)\sqrt{n} + \alpha} = \frac{\lambda}{\alpha + \lambda}.  \label{eq:boundlog2:tilde}
\end{eqnarray}
Combining \eqref{eq:expressionlog:tilde}, \eqref{eq:boundlog1:tilde}, and \eqref{eq:boundlog2:tilde} with the fact that $\frac{(\alpha+1) \tilde{\gamma} - \alpha n}{(\alpha+1)\tilde{\gamma} } >0$ gives
\begin{eqnarray*}
\text{\eqref{eq:expressionlog:tilde}}
&\geq& - \frac{1}{(1-\lambda)\alpha} \left( (\alpha+1)^2 + \frac{\lambda\alpha^2}{\alpha+ \lambda} \right)
 \cdot \int_{0}^{1} \frac{t}{1\cdot \left( 1- \frac{\lambda}{\alpha + \lambda} \right)} \, dt  \\
 &=& - \frac{1}{2} \,\frac{\alpha + \lambda}{(1-\lambda)\alpha^2}\left( (\alpha+1)^2 + \frac{\lambda\alpha^2}{\alpha + \lambda} \right),
\end{eqnarray*}
which implies
\begin{eqnarray}
\left(\frac{n}{\frac{\alpha+1}{\alpha}\tilde{\gamma}}\right)^{\tilde{\gamma}} \, \left(\frac{n}{(\alpha+1)(n-\tilde{\gamma})}\right)^{n-\tilde{\gamma}}
\geq \exp\left( - \frac{1}{2} \,\frac{\alpha + \lambda}{(1-\lambda)\alpha^2}\left( (\alpha+1)^2 + \frac{\lambda\alpha^2}{\alpha + \lambda} \right) \right). \label{eq:3.2}
\end{eqnarray}
Combining \eqref{eq:expansion:gamma:tilde}, \eqref{eq:3.1}, and \eqref{eq:3.2} completes the proof.
\end{proof}
\bigskip

\noindent \textbf{Proof of Proposition \ref{prop:upperboundforcoefficient}}

\begin{proof}
Let
\begin{equation*}
\beta = \left[ \frac{P^{\frac{1}{d}} (n+1) }{ P^{\frac{1}{d}} + 1 } \right].
\end{equation*}
As shown in Proposition \ref{proposition4},
\begin{equation*}
{n \choose k} \left( P^{\frac{1}{d}} \right)^{k} \leq {n \choose \beta} \left( P^{\frac{1}{d}} \right)^{\beta},
 \quad \text{for $k \in \{0, 1, \cdots, n \}$}.
\end{equation*}
Hence,
\begin{eqnarray}
\sum_{k=0}^{n} {n \choose k}^{d} P^{k}
&=& \sum_{k=0}^{n} \left({n \choose k} \left(P^{\frac{1}{d}}\right)^{k} \right)^{d-1} \,{n \choose k} \left(P^{\frac{1}{d}}\right)^{k} \notag \\[1mm]
&\leq& \sum_{k=0}^{n} \left({n \choose \beta} \left(P^{\frac{1}{d}}\right)^{\beta} \right)^{d-1} \,{n \choose k} \left(P^{\frac{1}{d}}\right)^{k}  \notag \\ [1mm]
&=& \left({n \choose \beta} \left(P^{\frac{1}{d}}\right)^{\beta} \right)^{d-1} \sum_{k=0}^{n} {n \choose k} \left(P^{\frac{1}{d}}\right)^{k} \notag \\ [1mm]
&=& \left({n \choose \beta} \left(P^{\frac{1}{d}}\right)^{\beta} \right)^{d-1} \left(1+P^{\frac{1}{d}}\right)^{n}. \label{eq:4.1}
\end{eqnarray}
By Proposition \ref{prop:upperboundforalpha}, we have
\begin{equation*}
{n \choose \beta} \leq \frac{M\left( P^{\frac{1}{d}} , \lambda \right)}{\sqrt{2 \, \pi}}\,n^{-\frac{1}{2}} \left(1 + P^{\frac{1}{d}}\right)^{n}.
\end{equation*}
The inequality in \eqref{eq:upperboundforcoefficients} now follows by plugging the above result into \eqref{eq:4.1}.
\end{proof}
\bigskip

\noindent \textbf{Proof of Proposition \ref{prop:lowerboundforcoefficient}}
\begin{proof}
Set
\begin{eqnarray*}
\beta = \left[ \frac{P^{\frac{1}{d}} (n+1) }{ P^{\frac{1}{d}} + 1 } \right], \
\gamma = \left[ \frac{P^{\frac{1}{d}} (n+1) }{ P^{\frac{1}{d}} + 1 } - \sqrt{n} \right], \  \text{and }
\tilde{\gamma} = \left[ \frac{P^{\frac{1}{d}} (n+1) }{ P^{\frac{1}{d}} + 1 } + \sqrt{n} \right].
\end{eqnarray*}
As shown in Proposition \ref{proposition4}, we have
\begin{eqnarray*}
&& {n \choose k} \left(P^{\frac{1}{d}}\right)^{k} \geq {n \choose \gamma} \left(P^{\frac{1}{d}}\right)^{\gamma}, \quad k \in \{ \gamma, \gamma +1, \cdots, \beta-1 \},  \\
&& {n \choose k} \left(P^{\frac{1}{d}}\right)^{k} \geq {n \choose \tilde{\gamma}} \left(P^{\frac{1}{d}}\right)^{\tilde{\gamma}}, \quad k \in \{ \beta, \beta +1, \cdots, \tilde{\gamma} \}.
\end{eqnarray*}
Hence,
\begin{eqnarray}
&& \sum_{k=0}^{n} {n \choose k}^{d} P^{k}
  = \sum_{k=0}^{n} \left({n \choose k} \left(P^{\frac{1}{d}} \right)^{k} \right)^{d} \notag \\ [1mm]
&\geq& \sum_{k=\gamma}^{\beta-1} \left({n \choose k} \left(P^{\frac{1}{d}} \right)^{k} \right)^{d}
       \, + \, \sum_{k=\beta}^{\tilde{\gamma}} \left({n \choose k} \left(P^{\frac{1}{d}} \right)^{k} \right)^{d}\notag \\ [1mm]
&\geq& \sum_{k=\gamma}^{\beta-1} \left({n \choose \gamma} \left(P^{\frac{1}{d}} \right)^{\gamma} \right)^{d}
       \, + \, \sum_{k=\beta}^{\tilde{\gamma}} \left({n \choose \tilde{\gamma}} \left(P^{\frac{1}{d}} \right)^{\tilde{\gamma}} \right)^{d}
\notag \\ [1mm]
&=& (\beta-\gamma)  \left({n \choose \gamma} \left(P^{\frac{1}{d}} \right)^{\gamma} \right)^{d} \,+\,
 (\tilde{\gamma} - \beta +1)  \left({n \choose \tilde{\gamma}} \left(P^{\frac{1}{d}} \right)^{\tilde{\gamma}} \right)^{d}.  \label{eq:5.1}
\end{eqnarray}
For $n \geq L\left(P^{\frac{1}{d}}, \lambda\right) $, the assumption in Proposition \ref{prop:lowerboundyamma} is satisfied when $\alpha$ is replaced by $P^{\frac{1}{d}}$. Applying Proposition \ref{prop:lowerboundyamma} and replacing $\alpha$ by $P^{\frac{1}{d}}$, we obtain
\begin{equation}
{ n \choose \gamma} \left( P^{\frac{1}{d}} \right)^{\gamma}
\geq \frac{\sqrt{2 \, \pi}}{e^2} \, C_{1}\left( P^{\frac{1}{d}} , \lambda \right) \, n^{-\frac{1}{2}} \left( 1 + P^{\frac{1}{d}} \right)^{n}.  \label{eq:5.2}
\end{equation}
Similarly, by Proposition \ref{prop:lowerboundyammatilde}, we have
\begin{equation}
{ n \choose \tilde{\gamma}} \left( P^{\frac{1}{d}} \right)^{\tilde{\gamma}}
\geq \frac{\sqrt{2 \, \pi}}{e^2} \, C_{2}\left( P^{\frac{1}{d}} , \lambda \right) \, n^{-\frac{1}{2}} \left( 1 + P^{\frac{1}{d}} \right)^{n}.  \label{eq:5.3}
\end{equation}
The condition $n \geq L\left(P^{\frac{1}{d}}, \lambda\right)$ implies
\begin{eqnarray}
 \sqrt{n} \geq \frac{ \left(P^{\frac{1}{d}}+1\right) + \sqrt{\left(P^{\frac{1}{d}}+1\right)^2 + 4\lambda P^{\frac{1}{d}}} }{2\lambda P^{\frac{1}{d}}}
 > \frac{P^{\frac{1}{d}}+1}{\lambda P^{\frac{1}{d}}}. \label{eq:5.3.5}
\end{eqnarray}
Then
\begin{eqnarray}
&&\beta - \gamma = \left[ \frac{P^{\frac{1}{d}} (n+1) }{ P^{\frac{1}{d}} + 1 } \right] - \left[ \frac{P^{\frac{1}{d}} (n+1) }{ P^{\frac{1}{d}} + 1 } -\sqrt{n} \right]  \notag \\  [1mm]
&>& \frac{P^{\frac{1}{d}} (n+1) }{ P^{\frac{1}{d}} + 1 } -1 - \left( \frac{P^{\frac{1}{d}} (n+1) }{ P^{\frac{1}{d}} + 1 } -\sqrt{n} \right) = \sqrt{n} -1  \notag  \\  [1mm]
&=& \sqrt{n} \cdot \left( 1- \frac{1}{\sqrt{n}} \right)
> \frac{(1-\lambda) P^{\frac{1}{d}} + 1}{P^{\frac{1}{d}} + 1} \, \sqrt{n}. \qquad (\text{applying \eqref{eq:5.3.5}})
\label{eq:5.4}
\end{eqnarray}
And
\begin{eqnarray}
&& \tilde{\gamma} - \beta + 1 = \left[ \frac{P^{\frac{1}{d}} (n+1) }{ P^{\frac{1}{d}} + 1 } + \sqrt{n} \right] - \left[ \frac{P^{\frac{1}{d}} (n+1) }{ P^{\frac{1}{d}} + 1 } \right] + 1  \notag \\  [1mm]
&>& \left( \frac{P^{\frac{1}{d}} (n+1) }{ P^{\frac{1}{d}} + 1 }  + \sqrt{n} - 1 \right) -\frac{P^{\frac{1}{d}} (n+1) }{ P^{\frac{1}{d}} + 1 } + 1 = \sqrt{n}.  \label{eq:5.5}
\end{eqnarray}
Combining \eqref{eq:5.1}, \eqref{eq:5.2}, \eqref{eq:5.3}, \eqref{eq:5.4}, and \eqref{eq:5.5} completes the proof.
\end{proof}

\end{document}